\documentclass{article}
\usepackage{amsmath, amssymb}
\usepackage{amssymb,pb-diagram,epsfig}
\usepackage{amscd}
\usepackage{array}
\usepackage{fancybox}

\newcommand{\ZZ}{\mathbb Z}
\newcommand{\PP}{\mathbb P}
\newcommand{\QQ}{\mathbb Q}
\newcommand{\CC}{\mathbb C}

\newcommand{\mcQ}{{\mathcal Q}}
\newcommand{\mcC}{{\mathcal C}}
\newcommand{\mcD}{{\mathcal D}}
\newcommand{\mcL}{{\mathcal L}}

\newcommand{\mcB}{{\mathcal B}}

\newcommand{\MW}{\mathop {\rm MW}\nolimits}

\newcommand{\rank}{\mathop {\rm rank}\nolimits}

\newcommand{\Sing}{\mathop {\rm Sing}\nolimits}

\newcommand{\I}{\mathop {\rm I}\nolimits}
\newcommand{\II}{\mathop {\rm II}\nolimits}
\newcommand{\III}{\mathop {\rm III}\nolimits}
\newcommand{\IV}{\mathop {\rm IV}\nolimits}

\newtheorem{thm}{Theorem}[section]
\newtheorem{mthm}{Theorem}

\newtheorem{cor}{Corollary}[section]
\newtheorem{prop}{Proposition}[section]

\newtheorem{lem}{Lemma}[section]
\newtheorem{defin}{Definition}[section]
\newtheorem{exmple}{Example}[section]
\newtheorem{rem}{Remark}[section]
\newtheorem{qz}{Question}[section]

\newtheorem{prbm}{Problem}[section]

\renewcommand{\thesubparagraph}{\theparagraph.\@arabic\c@subparagraph}
\begin{document}
 \begin{center}
  
 {\bf  \Large 
 Geometry of weak contact conics to irreducible quartics with 2 nodes and 1 cusp via rational elliptic surfaces and Zariski pairs
}
\bigskip

\bigskip
\large 
 Khulan TUMENBAYAR 
\end{center}
\bigskip

\small
{\bf Abstract.} Let $\mcQ$ be an irreducible quartic with two nodes and one cusp as its singularities and let $\mcC$ be a conic such that the intersection multiplicity at each point of $\mcC \cap \mcQ$ is even and $\mcC \cap \mcQ$ contain at least one smooth point $z_o$ of $\mcQ$. In this paper, for every $\mcQ$ we find all possible conics $\mcC$ as above via studying geometry of $\mcC$ and $\mcQ$ through that of integral sections of a rational elliptic surface which canonically arises from $\mcQ$ and $z_o \in \mcC \cap \mcQ$. As an application, we construct Zariski pairs of degree 7 and degree 8, whose irreducible components consist of $\mcQ$, $\mcC$ and line passing through two of the singular points of $\mcQ$ .



{\bf proposed running head:} Weak contact conics to certain quartics and Zariski pairs 
\normalsize
\section{Introduction} 
Let $\mcB$ be a reduced plane curve and let $\mcC$ be a smooth conic in $\PP^2$. Let $I_x(\mcB,\mcC)$ denote the intersection multiplicity at $x \in \mcB \cap \mcC$. We  first define a weak contact conic as follows:
\begin{defin}
 If $I_x(\mcB,\mcC) $ is even for $\forall x \in \mcB \cap \mcC$, then $\mcC$ is called a weak contact conic to $\mcB$. 
 Moreover if $\forall x \in \mcB \cap \mcC$ is a smooth point of $\mcB$ and $I_x(\mcB,\mcC) $ is even, then $\mcC$ is called a  contact conic to $\mcB$. 
\end{defin}
 
In \cite{khulan-tokunaga} and \cite{khulan-tokunaga2}, we studied contact conics to 3-nodal and 3-cuspidal quartics via geometry of rational elliptic surfaces. In this article, we continue to study geometry of contact conics to an irreducible quartic $\mcQ$ along this line. More precisely, we study weak contact conics to an irreducible quartic $\mcQ$ with 2 nodes and 1 cusp only as its singularities. 
In this article, we restrict ourselves to the case when $\mcC$ is a weak contact conic to $\mcQ$ such that $\mcQ \cap \mcC$ contains at least one smooth point $z_o$ of $\mcQ$ . Even under such restriction, we still obtain non-trivial examples of Zariski pairs, which we explain later. We first consider the following problem:
\begin{prbm}
Find all weak contact conics  that are tangent to $\mcQ$ at $z_o$.
\end{prbm}

 Let $x_1, x_2$ and $y$ denote the 2 nodes and the cusp of $\mcQ$, respectively. Let $(\mcQ \cap \mcC)_{sing}$ be the set of singular points of $\mcQ$ that are also in $\mcC$. For a weak contact conic passing through a smooth point $z_0$, we call it of type $i$, $i=1\dots 6$, if $(\mcQ \cap \mcC)_{sing}$ is as follows:
\begin{center}
\begin{tabular}{|c|c|c|c|c|c|c|}\hline
& type 1 & type 2& type 3 & type 4 & type 5 & type 6\\ \hline
$(\mcQ \cap \mcC)_{sing}$ & $\{y\}$ & $\{x_1\}$ or $\{x_2\}$ & $\{x_1,y\}$ or $\{x_2,y\}$ & $\{x_1,x_2\}$ &  $\{x_1,x_2,y\}$ & $\emptyset$   \\ \hline
\end{tabular} 
\end{center}
 
Let $l_{z_o}$ be the tangent line at $z_o$ to $\mcQ$. To answer to Problem 1.1, we consider 4 cases depending on the position of $l_{z_o}$ as follows:

\medskip
\begin{itemize}
\item $s$: $I_{z_o}(l_{z_o}, \mcQ) = 2$ or $3$, and $l_{z_o}$ meets $\mcQ$ transversely at other point(s).
\item $b$: $l_{z_o}$ is either bitangent line through $z_o$ or $I_{z_o}(l_{z_o}, Q) = 4$.
\item $sc$: $I_{z_o}(l_{z_o}, \mcQ) = 2$ and $l_{z_o}$ passes through the cusp of $\mcQ$.
\item $sn$: $I_{z_o}(l_{z_o}, \mcQ) = 2$ and $l_{z_o}$ passes through one of 2 nodes of $\mcQ$.
\end{itemize}

\begin{mthm} Let  $n_i$ be the number of weak contact conics of type $i$, $i=1\dots 6$. Then we have the following table:
\medskip
\begin{center}
\begin{tabular}{|c|c|c|c|c|c|c|c|} \hline 
 & $l_{z_o}\cap\mcQ$ & $n_1$ & $n_2$ & $n_3$ & $n_4$ & $n_5$& $n_6$\\ \hline
{\rm (I)} & s & 3 & 4 & 4 & 1 & 1&1  \\  
{\rm (II)} & b & 1 & 2 & 2 & 0 & 1 & 0 \\  
{\rm (III)} & sc & 0 & 2 & 0 & 1 & 0 & 0 \\ 
{\rm (IV)} & sn & 1 & 0 & 2 & 0 & 0 & 1\\ \hline
\end{tabular}
\end{center}                                                                                             
\end{mthm}

For $(\mcQ,z_o)$, likewise in \cite{bannai-tokunaga}, \cite{tokunaga10}, \cite{tokunaga14} and \cite{khulan-tokunaga2}, we consider a rational elliptic surface $S_{\mcQ,z_o}$ associated to a reduced plane quartic $\mcQ$ and its smooth point $z_o$.
Let $\MW(S_{\mcQ,z_o})$ be the set of sections of $\varphi : S_{\mcQ,z_o} \to \PP^1$. Our problem to find weak contact conics can  be reduced to study weak contact conics in $\MW(S_{\mcQ,z_o})$. As an application, we give examples of Zariski pairs, whose definition is as follows: 

\begin{defin}
A pair $(\mcB_1, \mcB_2)$ of reduced plane curves $\mcB_i~(i=1,2)$ of degree $n$ in $\PP^2$ is called a Zarski pair of degree $n$ if it satisfies the following condition:
\begin{itemize}
\item[(i)] $\mcB_i~(i= 1,2)$ are curves of degree $n$ such that the combinatorial type of $\mcB_1$ is the same as that of $\mcB_2$.
\item[(ii)] $(\PP^2,\mcB_1)$ is not homeomorphic to $(\PP^2,\mcB_2)$.
\end{itemize}

\end{defin}
For details, see \cite{act}. As for terminologies and notation concerning Zariski pairs, we use those used in \cite{act}, freely.
Now let us consider arrangements for Zariski pairs of degree 7 and 8, where irreducible components consist of quartic, conic and line as follows: 
\begin{itemize}
\item Let $\mcQ$ be an irreducible quartic with 2 nodes $x_1,~x_2$ and cusp $y$. 
\item Let $\mcL_i$ be line passing through $x_i$ and $y$, $i=1,2$. 
\item Let $\mcC_i$ be weak contact conic of type 1 for case $\rm (I)$, $i=0,1,2$.
\item Let $\bar{\mcC}$ be weak contact conic of type 5 for case $\rm (I)$.
\end{itemize}
\medskip
\textbf{Line-conic-quartic arrangements}

We put $\mcB_{i,j}:=\mcQ+\mcL_i+\mcC_j$, $i=1,2$, $j=0,1,2$.
\begin{enumerate}
\item[\textbf{1.}] 

 $\mcB_{1,j}$ and $\mcB_{2,j}$ have the same combinatorics, $(j=0,1,2)$.
 
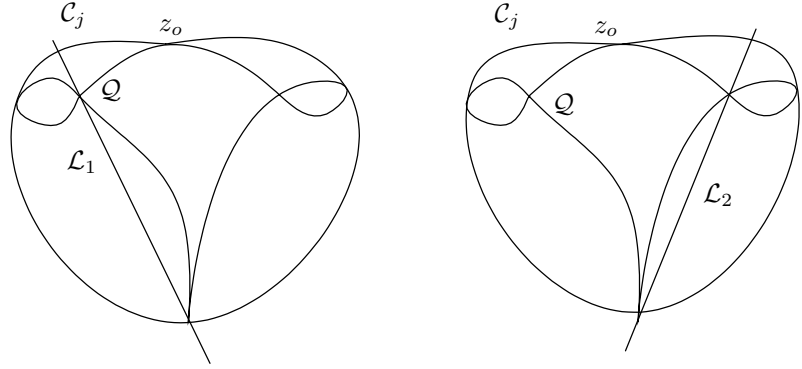
\begin{center}

{\unitlength 0.1in%
} 
 \end{center}
 
\item[\textbf{2.}]  For $(j\neq k, j=0,1,2, k=0,1,2)$, $\mcB_{i,j}$ and $\mcB_{i,k}$ have the same combinatorics, if 
$\mcQ \cap \mcC_j$ and $\mcQ \cap \mcC_k$ have the same number of points and there exists bijection $\varphi_{\rm{Sing}}$ between $\mcQ\cap\mcC_j$ and $\mcQ \cap \mcC_k$ such that for $\forall x \in \mcQ\cap \mcC_j$, $x$ and $\varphi_{\rm{Sing}}(x)$ have the same local topological types. 

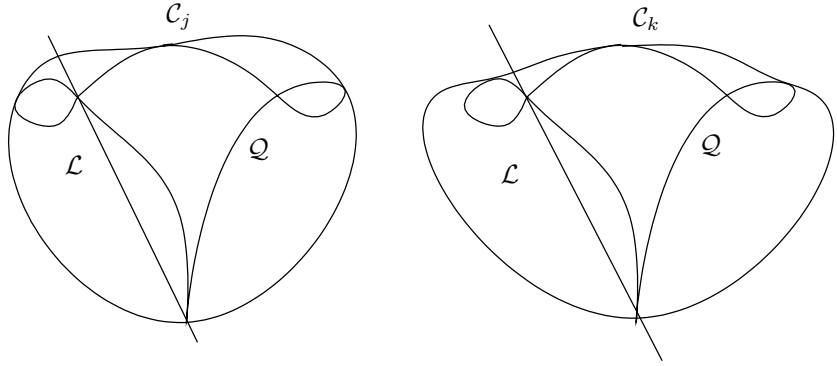
\begin{center}

{\unitlength 0.1in%
%
}
 \end{center}

\end{enumerate}

\textbf{Conic-conic-quartic arrangement}

We put $\mcD_{j}:=\mcQ+\bar{\mcC}+\mcC_j$, $j=0,1,2$. For $j\neq k, (j=0,1,2, k=0,1,2)$, $\mcD_{j}$ and $\mcD_{k}$ have the same combinatorics, if $\mcQ \cap \mcC_j$ and $\mcQ \cap \mcC_k$ have the same number of points and there exists bijection $\varphi_{\rm{Sing}}$ between $\mcQ\cap\mcC_j$ and $\mcQ \cap \mcC_k$ such that for $\forall x \in \mcQ\cap \mcC_j$, $x$ and $\varphi_{\rm{Sing}}(x)$ have the same local topological types. 

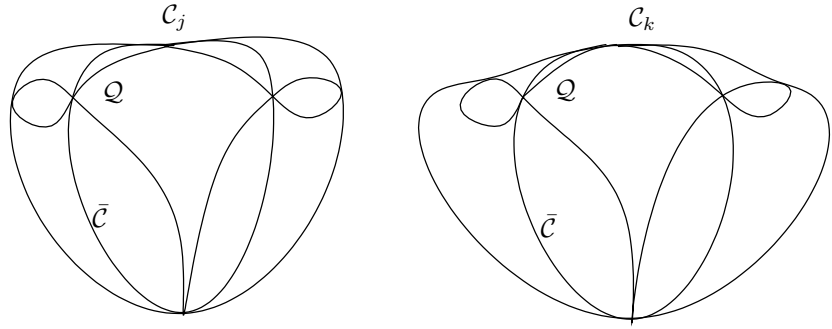
\begin{center}
 
{\unitlength 0.1in%
%
}%
 \end{center}

\begin{mthm}
Let $\mcB_{i,j}$ and $\mcD_k$ be arrangements defined above, $(i=1,2,~j=0,1,2,~k=0,1,2)$.
If 3 weak contact conics of type 1 for case $\rm (I)$ are chosen suitably by $\mcC_0,~\mcC_1,~\mcC_2$, then 
\begin{enumerate}
\item[(i)] $(\mcB_{1,1},\mcB_{2,1})$ and $(\mcB_{2,2},\mcB_{1,2})$ are Zariski pairs.  
\item[(ii)] If any pair of $(\mcB_{1,1},\mcB_{1,2})$, $(\mcB_{1,1},\mcB_{1,0})$, $(\mcB_{2,2},\mcB_{2,1})$, $(\mcB_{2,2},\mcB_{2,0})$, $(\mcD_0,\mcD_1)$ and $(\mcD_0,\mcD_2)$ has the same combinatorics then it is a Zariski pair.
\end{enumerate}
\end{mthm}

We proof this theorem in Section 4 and will give explicit examples of Zariski pairs in Section 5.

\section{Preliminaries}
\subsection{Construction of irreducible quartics with 2 nodes and 1 cusp}

 Let $[T,X,Z]$ be homogeneous coordinates of $\mathbb{P}^2$. 
 Let $Q$ be the standard quadratic transformation or the standard Cremona
 transformation with respect to $\{T=0\}, \{X=0\} $ and $\{Z=0\}$.
 We call $[0,0,1],[0,1,0]$ and $[1,0,0]$, the fundamental points with
 respect to $Q$.
 
\begin{lem}\label{lem:construction}
\begin{enumerate}
\item[(i)] Let C be a conic not passing through any of the fundamental points and tangent to $\{Z=0\}$, but not tangent to $\{T=0\}$ and $\{X=0\} $ in $\mathbb{P}^2$.
 Then $Q(C)$ is a quartic whose singularities are only 1 cusp at $[0,0,1]$ and 2 nodes at $[0,1,0]$ and $[1,0,0]$.

\item[(ii)] Let $L$ be the line tangent to $C$ at a point $P=[T_0,X_0,Z_0]\in C$, where $T_0X_0Z_0\neq 0$. If $L$ does not 
  contain any of the fundamental points, then $Q(L)$ is a conic tangent to $Q(C)$ at
 $Q(P)=[X_0Z_0,T_0Z_0,T_0X_0]$ and passes through the fundamental points.
 
\item[(iii)]Let $L$ be the line tangent to $C$ at a point $P=[T_0,X_0,Z_0]\in C$, where $T_0X_0Z_0\neq 0$.
    If $L$ pass through one of the fundamental points, then $Q(L)$ is a line passing through that fundamental point.
 
\item[(iv)] Let $L$ be conic, that contains  the fundamental points, then $Q(L)$ is a line.
 
\item[(v)] If $x\in \mathbb{P}^2 \setminus \{fundamental~points\}$,
  then $I_x(C,L)=I_{Q(x)}(Q(C),Q(L))$.
 
\end{enumerate}

\end{lem}

Since all of these statements are well-known, we omit their
proofs. We make use of Lemma~\ref{lem:construction} when
we consider explicit examples. 

\subsection{Elliptic Surfaces}

Throughout this article, an elliptic surface always means a smooth projective
surface $S$ with a fibration $\varphi : S \to C$ over a smooth projective
curve, $C$, as follows:
\begin{enumerate}

\item[(i)] There exists
non empty finite subset  $\Sing(\varphi) \subset C$ such that
$\varphi^{-1}(v)$ is a smooth curve of genus $1$ for $v \in
C\smallsetminus \Sing(\varphi)$, while $\varphi^{-1}(v)$ is not a
smooth curve of genus $1$ for $v \in \Sing(\varphi)$.

\item[(ii)] There exists a
section $O: C \to S$ (we identify $O$ with its image in $S$).

\item[(iii)]
there is no exceptional curve of the first kind in any fiber.
\end{enumerate}

For $v \in \Sing(\varphi)$, we call $F_v = \varphi^{-1}(v)$ a singular fiber over $v$. As
for  the types of singular fibers, we
use notation given by
Kodaira (\cite{kodaira}).
We denote the irreducible  decomposition of $F_v$ by
\[
F_v = \Theta_{v, 0} + \sum_{i=1}^{m_v-1}a_{v,i}\Theta_{v,i},
\]
where $m_v$ is the number of irreducible components of $F_v$ and
$\Theta_{v,0}$ is the irreducible component with $\Theta_{v,0}O = 1$. 


Let $\MW(S)$ be the set of sections of $\varphi : S \to C$. By our assumption, $\MW(S) \neq \emptyset$.
It is known that $\MW(S)$ has the structure of an abelian group with zero element $O$. It is called the Mordell-Weil group.
We denote its addition and the multiplication-by-$m$ map by $s_1 \dot{+} s_2$ and $[m]s_1 $, respectively, where $s_1, s_2 \in \MW(S)$ $m \in \ZZ$.
 
Let $E$ be a generic fiber of elliptic surface $\varphi : S \to C$. One can regard $E$
as an elliptic curve over ${\mathbb C}(C)$, the rational function field of $C$.
It is known that $\MW(S)$ can be identified with the set of ${\mathbb C}(C)$-rational points $E(\CC(C))$.
We also denote the addition and the multiplication-by-$m$ map on $E(\CC(C))$ by $P_1 \dot{+} P_2$ and $[m]P_1$ for
$P_1, P_2 \in E(\CC(C))$, respectively.

In \cite{shioda90}, Shioda introduced a $\QQ$-valued bilinear form on $E(\CC(C))$ called the height pairing. We denote
it by $\langle \, , \, \rangle$. For an explicit formula of $\langle P_1,
P_2\rangle$ ($P_1, P_2 \in E(\CC(C))$), see \cite[Theorem 8.6]{shioda90}.

\subsection{Rational Elliptic Surfaces $S_{\mcQ,z_o}$ }
We refer to \cite{tokunaga14} for details. Let $\varphi:S \to \PP^1$ be a rational elliptic surface and let $\mathcal{E}$ be the generic fiber of $\varphi$. The inverse morphism on $\mathcal{E}$ with respect to the group law on $\mathcal{E}$ induces an involution $[-1]_{\varphi}$ on $S$. Let $S/ \langle[-1]_{\varphi}\rangle$ be the quotient by $[-1]_{\varphi}$. It is known that $S/ \langle[-1]_{\varphi}\rangle$ is smooth and if $\varphi$ is a rational elliptic surface and $S$ has a reducible singular fiber then we can blow down $S/ \langle[-1]_{\varphi}\rangle$ to $\PP^2$. 
Hence we have the following double cover diagram

\[
\begin{CD}
S'' @<{\bar{\mu}}<< S \\
@V{f''}VV                 @VV{f}V \\
\PP^2@<<\bar{q}< \widehat{\Sigma}_2,
\end{CD}
\]
where 
\begin{itemize}
\item $f:S \to S/ \langle[-1]_{\varphi}\rangle=\widehat{\Sigma}_2$: the quotient morphism,
\item $\bar{q}: \widehat{\Sigma}_2 \to \PP^2$, where $\bar{q}$ is a composition of blowing ups, 
\item $S\xrightarrow{\bar{\mu}} S'' \xrightarrow{f''} \PP^2$: the Stein factorization of $\bar{q}\circ f$.
\end{itemize}

 Note that $S''$ is a double cover of $\PP^2$ with branch locus $\mcQ$, a quartic, as explained in \cite{tokunaga14}.
 
\medskip
Conversely, let $\mcQ$ be a reduced quartic, not concurrent 4 lines and let $z_o$ be a smooth point on $\mcQ$. Let $f_{\mcQ,z_o}:S \to \PP^2$ be the composition of maps as follows:
\[
f_{\mcQ,z_o}:S \xrightarrow{\pi} S' \xrightarrow{f'} W \xrightarrow{\bar{q}_2} \widehat{\PP^2} \xrightarrow{\bar{q}_1} \PP^2
\]
where

\begin{itemize}
\item $\bar{q}_1:~\widehat{\PP^2} \to \PP^2$ blow-up at $z_o$ twice,
\item $\bar{q}_2:~W \to \widehat{\PP^2} $ blow-up until $\widehat{\mcQ}$, the strict transform of $\mcQ$, becomes smooth,
\item $f':~S' \to W$ the double cover of $W$ with branch locus $\widehat{\mcQ}+\Delta_0$, where $\Delta_0$ is the strict image of the exceptional curve of first blow up of $\bar{q}_1$,
\item $\pi:~S \to S'$ the minimal resolution of singularities of $S'$. 

\end{itemize}
If $L$ is line passing through $z_o$ which is not equal to $l_{z_o}$ then $f^{-1}_{\mcQ,z_o}(L)$ is curve of genus 1 on $S$. 
Hence we have a rational elliptic surface $\varphi : S_{\mcQ,z_o} \to \PP^1$.

\subsection{A criterion for Zariski pair}
From now on, let $\mcQ$ be an irreducible quartic which has 1 cusp and 2 nodes only as its singularities.
We fix an isomorphism $\MW(S)\cong M_0 \oplus \MW_{tor}$, $M_0 \cong \ZZ^{\oplus r}$, $r=\rank \MW(S)$.
Choose $s_1,s_2 \in M_0$ so that $s_1$ and $s_2$ are a part of a basis of $\ZZ^{\oplus r}$, i.e.,$M_0/\ZZ s_1+\ZZ s_2$ is free of rank $r-2$. Put $\bar{\mcC}_i:=f_{\mcQ,z_o}(s_i)$ and $\mcC_i:=f_{\mcQ,z_o}([2]s_i)$ $i=1,2$. Then we have the following theorem.
\begin{thm}[{\cite[Proposition 4.4]{tokunaga14}}]\label{thm:homeo} 
Let $\mcQ$, $\bar{\mcC}_i$ and $\mcC_i$, $(i=1,2)$ be as above. Then there is no homeomorphism $h: (\PP^2,\mcQ+\bar{\mcC}_1+\mcC_1) \to (\PP^2,\mcQ+\bar{\mcC_2}+\mcC_1)$, such that $h(\mcQ)=\mcQ$.
\end{thm}
Using the same idea as in \cite{tokunaga14}, we have similar theorem. More precisely, in Theorem \ref{thm:homeo} we used the fact that $s_1$ and $[2]s_1$ are linearly dependent, while $s_2$ and $[2]s_1$ are linearly independent. In the following theorem we used the fact that $s_1$ and $[2]s_1$ are linearly dependent, while $s_1$ and $s_2$ are linearly independent.
\begin{thm}\label{thm:homeo1} 
Let $\mcQ$, $\bar{\mcC}_i$ and $\mcC_i$, $(i=1,2)$ be as above. Then there is no homeomorphism $h: (\PP^2,\mcQ+\bar{\mcC}_1+\mcC_1) \to (\PP^2,\mcQ+\bar{\mcC}_1+\mcC_2)$, such that $h(\mcQ)=\mcQ$.
\end{thm}
\section{Computation of weak contact conics in $S_{\mcQ, z_o}$}

Let $\mcQ$ be the same as above and let $z_o$ be a smooth point on $\mcQ$.
The tangent line $l_{z_o}$ to $\mcQ$ at $z_o$ gives rise to a singular fiber of $\varphi$ whose type is determined by
how $l_{z_o}$ intersects with $\mcQ$ as follows:
\begin{center}
\begin{tabular}{|c|c|l|}  \hline
(i) & $\I_2$ & $l_{z_o}$ meets $\mcQ$ with two other distinct points. \\ \hline
(ii) & $\III$ & $l_{z_o}$ is a $3$-fold tangent point. \\ \hline
(iii) &$\I_3$ & $l_{z_o}$ is a bitangent line. \\ \hline
(iv) & $\IV$ & $l_{z_o}$ is a $4$-fold tangent point. \\ \hline
(v)  & $\I_5$ & $l_{z_o}$ passes through the cusp of $\mcQ$  \\ \hline
(vi) & $\I_4$ & $l_{z_o}$ passes through a node of $\mcQ$  \\ \hline
\end{tabular}
\end{center}
From \cite[Table 6.2]{miranda-persson} and the above table,  we have the following table for possible configurations of singular fibers of $S_{\mcQ, z_o}$:
\begin{center}
\begin{tabular}{|c|c|c|} \hline
   & Singular fibers & the position of $l_{z_o}$ \\ \hline
(a) & $\{2\I_2, \I_3, \I_2, 3\I_1\}, \{2\I_2, \I_3, \I_2, \I_1, \II\}$ & (i)  \\ \hline
(b) & $\{2\I_2, \I_3, \III, 2\I_1\}, \{2\I_2, \I_3 \III, \II\}$ & (ii) \\ \hline
(c) & $\{2\I_2, 2\I_3, 2\I_1\}, \{2\I_2, 2\I_3, \II\}$ & (iii) \\ \hline
(d) & $\{2\I_2, \I_3, \IV,  \I_1\}$  & (iv) \\ \hline
(e) & $\{2\I_2, \I_5, 3\I_1\}, \{2\I_2, \I_5, \I_1, \II\}$ & (v)  \\ \hline
(f) & $\{\I_2, \I_3, \I_4, 3\I_1\}, \{\I_2, \I_3, \I_4, \I_1, \II\}$ & (vi)  \\ \hline
\end{tabular}
\end{center}
In our later argument, we need to know the structure of $E_{\mcQ, z_o}(\CC(t))$. 
Since irreducible singular fibers and
the difference between $\III$ (resp. $\IV$) type and $\I_2$ (resp. $\I_3$) type do not affect the
structure of $E_{\mcQ, z_o}(\CC(t))$, we only consider the cases:
\[
\mathrm{(I)}\, 2\I_2, \I_3, \I_2, 3\I_1, \quad
\mathrm{(II)}\, 2\I_2, 2\I_3, 2\I_1, \quad
\mathrm{(III)}\, 2\I_2, \I_5, 3\I_1, \quad
\mathrm{(IV)}\, \I_2, \I_3, \I_4, 3\I_1.
\]
In \cite{oguiso-shioda}, the structure of $E_{\mcQ, z_o}(\CC(t))$ with respect to the height pairing $\langle~,~\rangle$ is given by Gram-matrices as follows:
\[
\mathrm{(I)}\, A_1^*\oplus \frac 16 \left [\begin{array}{cc}
                                                                2 & 1 \\
                                                                1 & 2
                                                                \end{array} \right ], \quad
\mathrm{(II)}\, \langle  1/6 \rangle^{\oplus 2},  \quad
\mathrm{(III)}\, \frac{1}{10} \left [\begin{array}{cc}
                                                                2 & 1 \\
                                                                1 & 3
                                                                \end{array} \right ], \quad
\mathrm{(IV)}\, A_1^* \oplus \langle 1/12 \rangle.
\]
Let $\mcL_1$ and $\mcL_2$ be two lines passing through the cusp and one of the two nodes of $\mcQ$ and 
$\mcL_3$ be line passing through two nodes of $\mcQ$. Let $\mcC_0$ be conic that passes through 3 singular points of $\mcQ$ and tangent to $\mcQ$ at $z_0$. For the cases $\mathrm{(I)},\mathrm{(II)}$ there exists such conic by lemma~\ref{lem:construction}(ii). Note that for the cases $\mathrm{(III)},\mathrm{(IV)}$ there exists no such conic. More precisely, if such $\mcC_0$ exists then $Q(l_{z_0})$ and $Q(\mcC_0)$ would be same line, which is contradiction.

By our construction of $S_{\mcQ, z_o}$, $\mcL_i$ ($i = 1, 2, 3$) and $\mcC_0$ give rise to 
sections $s_i^\pm$ ($i = 1, 2, 3$) and $s_0^\pm$, respectively. We also denote the corresponding element in
$E_{\mcQ, z_o}(\CC(t))$ by $P_i^\pm$ $(i=1,2,3)$ and $P_0^\pm$, respectively.

\medskip
{\bf Case $\mathrm{(I)}$}
 
 We label irreducible components of singular fibers of type $\I_2$ in such a way that
 $\Theta_{i, 1}$ ($i = 1, 2$) are those arising from the nodes of $\mcQ$ and  that of type $\I_3$ such that $\Theta_{3,1}$, $\Theta_{3,2}$ are irreducible components from the cusp of $\mcQ$. Also $\Theta_{\infty, 1}$
 is the one from $l_{z_o}$. 
 By our construction of $S_{\mcQ, z_o}$, we may assume that
 $s_i^+$ ($i = 0, 1, 2, 3$) meet each singular fiber as in the
 figure below.
  \begin{center}
{\unitlength 0.1in%
\begin{picture}( 45.4200, 21.1100)(  6.7500,-22.6600)%
%
\special{pn 13}%
\special{pa 1285 253}%
\special{pa 1256 273}%
\special{pa 1172 333}%
\special{pa 1145 353}%
\special{pa 1119 374}%
\special{pa 1093 396}%
\special{pa 1045 440}%
\special{pa 1023 463}%
\special{pa 1002 487}%
\special{pa 982 512}%
\special{pa 965 537}%
\special{pa 949 564}%
\special{pa 935 592}%
\special{pa 923 621}%
\special{pa 913 651}%
\special{pa 906 682}%
\special{pa 901 714}%
\special{pa 898 747}%
\special{pa 898 781}%
\special{pa 900 814}%
\special{pa 905 848}%
\special{pa 911 881}%
\special{pa 920 914}%
\special{pa 931 946}%
\special{pa 945 977}%
\special{pa 961 1006}%
\special{pa 979 1034}%
\special{pa 999 1060}%
\special{pa 1021 1083}%
\special{pa 1046 1105}%
\special{pa 1071 1125}%
\special{pa 1097 1145}%
\special{pa 1122 1164}%
\special{pa 1145 1185}%
\special{pa 1167 1206}%
\special{pa 1186 1229}%
\special{pa 1203 1254}%
\special{pa 1218 1279}%
\special{pa 1231 1306}%
\special{pa 1242 1334}%
\special{pa 1251 1363}%
\special{pa 1259 1393}%
\special{pa 1266 1425}%
\special{pa 1271 1456}%
\special{pa 1274 1489}%
\special{pa 1277 1523}%
\special{pa 1279 1557}%
\special{pa 1279 1627}%
\special{pa 1278 1662}%
\special{pa 1277 1698}%
\special{pa 1275 1734}%
\special{pa 1273 1771}%
\special{pa 1271 1807}%
\special{pa 1268 1844}%
\special{pa 1268 1846}%
\special{fp}%
%
\special{pn 13}%
\special{pa 1022 269}%
\special{pa 1044 296}%
\special{pa 1065 322}%
\special{pa 1087 349}%
\special{pa 1129 403}%
\special{pa 1169 457}%
\special{pa 1205 511}%
\special{pa 1223 539}%
\special{pa 1239 566}%
\special{pa 1254 594}%
\special{pa 1280 650}%
\special{pa 1291 678}%
\special{pa 1300 706}%
\special{pa 1308 735}%
\special{pa 1314 764}%
\special{pa 1318 793}%
\special{pa 1320 823}%
\special{pa 1321 852}%
\special{pa 1320 882}%
\special{pa 1318 912}%
\special{pa 1314 943}%
\special{pa 1310 973}%
\special{pa 1296 1035}%
\special{pa 1288 1066}%
\special{pa 1268 1128}%
\special{pa 1257 1159}%
\special{pa 1245 1191}%
\special{pa 1232 1222}%
\special{pa 1219 1254}%
\special{pa 1205 1286}%
\special{pa 1191 1317}%
\special{pa 1177 1349}%
\special{pa 1132 1445}%
\special{pa 1126 1457}%
\special{fp}%
%
\special{pn 13}%
\special{pa 1988 253}%
\special{pa 1959 273}%
\special{pa 1875 333}%
\special{pa 1848 353}%
\special{pa 1822 374}%
\special{pa 1796 396}%
\special{pa 1772 417}%
\special{pa 1748 440}%
\special{pa 1726 463}%
\special{pa 1705 487}%
\special{pa 1686 512}%
\special{pa 1668 537}%
\special{pa 1652 564}%
\special{pa 1638 592}%
\special{pa 1626 621}%
\special{pa 1616 651}%
\special{pa 1609 682}%
\special{pa 1604 714}%
\special{pa 1601 747}%
\special{pa 1601 781}%
\special{pa 1603 814}%
\special{pa 1607 848}%
\special{pa 1614 881}%
\special{pa 1623 914}%
\special{pa 1634 946}%
\special{pa 1648 977}%
\special{pa 1664 1006}%
\special{pa 1682 1034}%
\special{pa 1702 1060}%
\special{pa 1724 1083}%
\special{pa 1749 1105}%
\special{pa 1774 1125}%
\special{pa 1800 1144}%
\special{pa 1825 1164}%
\special{pa 1849 1184}%
\special{pa 1870 1206}%
\special{pa 1890 1228}%
\special{pa 1907 1253}%
\special{pa 1922 1278}%
\special{pa 1935 1305}%
\special{pa 1946 1333}%
\special{pa 1955 1362}%
\special{pa 1963 1392}%
\special{pa 1969 1424}%
\special{pa 1974 1455}%
\special{pa 1978 1488}%
\special{pa 1981 1522}%
\special{pa 1983 1590}%
\special{pa 1983 1626}%
\special{pa 1982 1661}%
\special{pa 1980 1697}%
\special{pa 1979 1734}%
\special{pa 1976 1770}%
\special{pa 1974 1807}%
\special{pa 1971 1843}%
\special{pa 1971 1846}%
\special{fp}%
%
\special{pn 13}%
\special{pa 1726 269}%
\special{pa 1748 296}%
\special{pa 1769 323}%
\special{pa 1791 349}%
\special{pa 1812 376}%
\special{pa 1832 403}%
\special{pa 1853 430}%
\special{pa 1891 484}%
\special{pa 1909 512}%
\special{pa 1926 539}%
\special{pa 1942 566}%
\special{pa 1957 594}%
\special{pa 1971 622}%
\special{pa 1983 650}%
\special{pa 1994 678}%
\special{pa 2003 707}%
\special{pa 2011 736}%
\special{pa 2017 765}%
\special{pa 2021 794}%
\special{pa 2023 823}%
\special{pa 2024 853}%
\special{pa 2023 883}%
\special{pa 2021 913}%
\special{pa 2017 943}%
\special{pa 2013 974}%
\special{pa 2006 1004}%
\special{pa 1999 1035}%
\special{pa 1991 1066}%
\special{pa 1981 1097}%
\special{pa 1971 1129}%
\special{pa 1960 1160}%
\special{pa 1948 1191}%
\special{pa 1922 1255}%
\special{pa 1909 1286}%
\special{pa 1894 1318}%
\special{pa 1880 1350}%
\special{pa 1835 1446}%
\special{pa 1830 1457}%
\special{fp}%
%
\special{pn 13}%
\special{pa 2387 253}%
\special{pa 2387 1851}%
\special{fp}%
%
\special{pn 13}%
\special{pa 3512 258}%
\special{pa 3484 278}%
\special{pa 3455 298}%
\special{pa 3399 338}%
\special{pa 3372 359}%
\special{pa 3346 379}%
\special{pa 3320 401}%
\special{pa 3272 445}%
\special{pa 3250 468}%
\special{pa 3229 492}%
\special{pa 3210 517}%
\special{pa 3192 543}%
\special{pa 3176 569}%
\special{pa 3162 597}%
\special{pa 3150 626}%
\special{pa 3141 656}%
\special{pa 3133 687}%
\special{pa 3128 719}%
\special{pa 3126 752}%
\special{pa 3125 786}%
\special{pa 3127 820}%
\special{pa 3132 853}%
\special{pa 3138 887}%
\special{pa 3147 919}%
\special{pa 3159 951}%
\special{pa 3172 982}%
\special{pa 3188 1011}%
\special{pa 3206 1039}%
\special{pa 3226 1065}%
\special{pa 3249 1089}%
\special{pa 3273 1110}%
\special{pa 3298 1130}%
\special{pa 3324 1150}%
\special{pa 3349 1169}%
\special{pa 3373 1190}%
\special{pa 3394 1211}%
\special{pa 3413 1234}%
\special{pa 3430 1259}%
\special{pa 3445 1284}%
\special{pa 3458 1311}%
\special{pa 3469 1339}%
\special{pa 3478 1369}%
\special{pa 3486 1399}%
\special{pa 3493 1430}%
\special{pa 3498 1462}%
\special{pa 3501 1494}%
\special{pa 3504 1528}%
\special{pa 3506 1562}%
\special{pa 3506 1632}%
\special{pa 3505 1667}%
\special{pa 3504 1703}%
\special{pa 3502 1740}%
\special{pa 3500 1776}%
\special{pa 3498 1813}%
\special{pa 3495 1849}%
\special{pa 3495 1851}%
\special{fp}%
%
\special{pn 13}%
\special{pa 3249 274}%
\special{pa 3271 301}%
\special{pa 3292 327}%
\special{pa 3314 354}%
\special{pa 3356 408}%
\special{pa 3396 462}%
\special{pa 3432 516}%
\special{pa 3450 544}%
\special{pa 3466 571}%
\special{pa 3481 599}%
\special{pa 3507 655}%
\special{pa 3518 683}%
\special{pa 3527 711}%
\special{pa 3535 740}%
\special{pa 3541 769}%
\special{pa 3545 798}%
\special{pa 3547 828}%
\special{pa 3548 857}%
\special{pa 3547 887}%
\special{pa 3545 917}%
\special{pa 3541 948}%
\special{pa 3537 978}%
\special{pa 3523 1040}%
\special{pa 3515 1071}%
\special{pa 3495 1133}%
\special{pa 3484 1164}%
\special{pa 3472 1196}%
\special{pa 3459 1227}%
\special{pa 3446 1259}%
\special{pa 3432 1291}%
\special{pa 3418 1322}%
\special{pa 3404 1354}%
\special{pa 3359 1450}%
\special{pa 3353 1462}%
\special{fp}%
%
\special{pn 13}%
\special{pa 4008 258}%
\special{pa 4014 290}%
\special{pa 4019 321}%
\special{pa 4031 385}%
\special{pa 4036 416}%
\special{pa 4046 480}%
\special{pa 4051 511}%
\special{pa 4059 575}%
\special{pa 4063 606}%
\special{pa 4069 670}%
\special{pa 4075 766}%
\special{pa 4075 862}%
\special{pa 4074 894}%
\special{pa 4066 1022}%
\special{pa 4063 1054}%
\special{pa 4061 1086}%
\special{pa 4052 1182}%
\special{pa 4050 1213}%
\special{pa 4047 1245}%
\special{pa 4043 1277}%
\special{pa 4028 1307}%
\special{pa 4002 1329}%
\special{pa 3971 1331}%
\special{pa 3938 1313}%
\special{pa 3911 1283}%
\special{pa 3897 1249}%
\special{pa 3901 1221}%
\special{pa 3926 1204}%
\special{pa 3962 1197}%
\special{pa 3996 1201}%
\special{pa 4023 1213}%
\special{pa 4044 1232}%
\special{pa 4059 1258}%
\special{pa 4069 1290}%
\special{pa 4076 1325}%
\special{pa 4078 1363}%
\special{pa 4079 1403}%
\special{pa 4078 1444}%
\special{pa 4076 1484}%
\special{pa 4074 1523}%
\special{pa 4072 1560}%
\special{pa 4070 1595}%
\special{pa 4069 1629}%
\special{pa 4068 1662}%
\special{pa 4066 1724}%
\special{pa 4064 1782}%
\special{pa 4064 1811}%
\special{pa 4063 1840}%
\special{pa 4063 1851}%
\special{fp}%
%
\special{pn 13}%
\special{pa 4397 264}%
\special{pa 4403 296}%
\special{pa 4408 327}%
\special{pa 4413 359}%
\special{pa 4419 391}%
\special{pa 4424 422}%
\special{pa 4434 486}%
\special{pa 4438 517}%
\special{pa 4443 549}%
\special{pa 4447 581}%
\special{pa 4450 613}%
\special{pa 4453 644}%
\special{pa 4456 676}%
\special{pa 4460 740}%
\special{pa 4462 804}%
\special{pa 4462 868}%
\special{pa 4460 932}%
\special{pa 4456 996}%
\special{pa 4454 1029}%
\special{pa 4451 1061}%
\special{pa 4449 1092}%
\special{pa 4440 1188}%
\special{pa 4438 1220}%
\special{pa 4435 1251}%
\special{pa 4430 1283}%
\special{pa 4415 1313}%
\special{pa 4389 1335}%
\special{pa 4357 1336}%
\special{pa 4325 1317}%
\special{pa 4299 1286}%
\special{pa 4285 1253}%
\special{pa 4291 1225}%
\special{pa 4317 1208}%
\special{pa 4353 1202}%
\special{pa 4387 1206}%
\special{pa 4414 1219}%
\special{pa 4434 1239}%
\special{pa 4448 1265}%
\special{pa 4458 1297}%
\special{pa 4464 1333}%
\special{pa 4466 1371}%
\special{pa 4466 1411}%
\special{pa 4465 1452}%
\special{pa 4462 1492}%
\special{pa 4460 1531}%
\special{pa 4459 1568}%
\special{pa 4457 1603}%
\special{pa 4456 1637}%
\special{pa 4454 1701}%
\special{pa 4454 1731}%
\special{pa 4453 1761}%
\special{pa 4453 1790}%
\special{pa 4452 1819}%
\special{pa 4452 1857}%
\special{fp}%
%
\special{pn 13}%
\special{pa 4762 258}%
\special{pa 4768 290}%
\special{pa 4773 321}%
\special{pa 4779 353}%
\special{pa 4784 385}%
\special{pa 4789 416}%
\special{pa 4799 480}%
\special{pa 4804 511}%
\special{pa 4812 575}%
\special{pa 4815 607}%
\special{pa 4819 638}%
\special{pa 4821 670}%
\special{pa 4824 702}%
\special{pa 4825 734}%
\special{pa 4827 766}%
\special{pa 4827 862}%
\special{pa 4826 894}%
\special{pa 4824 926}%
\special{pa 4823 958}%
\special{pa 4821 990}%
\special{pa 4818 1022}%
\special{pa 4816 1054}%
\special{pa 4813 1086}%
\special{pa 4811 1118}%
\special{pa 4805 1182}%
\special{pa 4803 1214}%
\special{pa 4800 1245}%
\special{pa 4796 1277}%
\special{pa 4780 1308}%
\special{pa 4755 1329}%
\special{pa 4723 1331}%
\special{pa 4691 1313}%
\special{pa 4664 1282}%
\special{pa 4651 1249}%
\special{pa 4656 1220}%
\special{pa 4681 1203}%
\special{pa 4717 1197}%
\special{pa 4751 1201}%
\special{pa 4778 1213}%
\special{pa 4798 1233}%
\special{pa 4813 1260}%
\special{pa 4823 1291}%
\special{pa 4828 1326}%
\special{pa 4831 1365}%
\special{pa 4831 1405}%
\special{pa 4830 1446}%
\special{pa 4827 1486}%
\special{pa 4825 1525}%
\special{pa 4824 1562}%
\special{pa 4822 1597}%
\special{pa 4821 1631}%
\special{pa 4819 1695}%
\special{pa 4819 1725}%
\special{pa 4818 1755}%
\special{pa 4818 1784}%
\special{pa 4817 1813}%
\special{pa 4817 1851}%
\special{fp}%
%
\special{pn 13}%
\special{pa 771 474}%
\special{pa 930 474}%
\special{fp}%
%
\special{pn 13}%
\special{pa 1072 474}%
\special{pa 1640 474}%
\special{fp}%
%
\special{pn 13}%
\special{pa 1787 480}%
\special{pa 2350 480}%
\special{fp}%
%
\special{pn 13}%
\special{pa 2404 480}%
\special{pa 3353 474}%
\special{fp}%
%
\special{pn 13}%
\special{pa 3462 480}%
\special{pa 4921 480}%
\special{fp}%
%
\special{pn 13}%
\special{pa 777 685}%
\special{pa 853 685}%
\special{fp}%
%
\special{pn 13}%
\special{pa 978 685}%
\special{pa 1928 685}%
\special{fp}%
%
\special{pn 13}%
\special{pa 2109 685}%
\special{pa 2338 679}%
\special{fp}%
%
\special{pn 13}%
\special{pa 2447 685}%
\special{pa 3070 685}%
\special{fp}%
%
\special{pn 13}%
\special{pa 3206 685}%
\special{pa 4926 685}%
\special{fp}%
%
\special{pn 13}%
\special{pa 771 868}%
\special{pa 1251 868}%
\special{fp}%
%
\special{pn 13}%
\special{pa 1377 863}%
\special{pa 1377 868}%
\special{fp}%
\special{pa 1372 868}%
\special{pa 1558 868}%
\special{fp}%
%
\special{pn 13}%
\special{pa 1661 868}%
\special{pa 2338 868}%
\special{fp}%
%
\special{pn 13}%
\special{pa 2288 1446}%
\special{pa 2769 712}%
\special{fp}%
%
\special{pn 13}%
\special{pa 2753 814}%
\special{pa 2333 274}%
\special{fp}%
%
\special{pn 13}%
\special{pa 2452 868}%
\special{pa 3087 868}%
\special{fp}%
%
\special{pn 13}%
\special{pa 3196 868}%
\special{pa 4926 868}%
\special{fp}%
%
\special{pn 13}%
\special{pa 777 1079}%
\special{pa 968 1079}%
\special{fp}%
%
\special{pn 13}%
\special{pa 1077 1079}%
\special{pa 1655 1079}%
\special{fp}%
%
\special{pn 13}%
\special{pa 1765 1079}%
\special{pa 2469 1079}%
\special{fp}%
%
\special{pn 13}%
\special{pa 2556 1079}%
\special{pa 3174 1079}%
\special{fp}%
%
\special{pn 13}%
\special{pa 3305 1079}%
\special{pa 4909 1074}%
\special{fp}%
%
\special{pn 13}%
\special{pa 782 1657}%
\special{pa 4926 1657}%
\special{fp}%
%
\special{pn 13}%
\special{pa 675 155}%
\special{pa 5217 155}%
\special{pa 5217 2266}%
\special{pa 675 2266}%
\special{pa 675 155}%
\special{pa 5217 155}%
\special{fp}%
\put(50.7000,-6.7000){\makebox(0,0){$s_1$}}%
\put(50.7000,-4.6500){\makebox(0,0){$s_0$}}%
\put(50.7500,-8.7000){\makebox(0,0){$s_2$}}%
\put(50.6000,-16.5500){\makebox(0,0){$O$}}%
\put(14.5000,-17.9000){\makebox(0,0){$\Theta_{1,0}$}}%
\put(21.5000,-18.0000){\makebox(0,0){$\Theta_{2,0}$}}%
\put(36.6000,-18.0500){\makebox(0,0){$\Theta_{\infty,0}$}}%
\put(14.0000,-11.9500){\makebox(0,0){$\Theta_{1,1}$}}%
\put(21.0000,-11.9500){\makebox(0,0){$\Theta_{2,1}$}}%
\put(37.2500,-11.9000){\makebox(0,0){$\Theta_{\infty,1}$}}%
\put(26.3000,-11.9000){\makebox(0,0){$\Theta_{3,2}$}}%
\put(25.6000,-18.0500){\makebox(0,0){$\Theta_{3,0}$}}%
\put(26.2000,-3.5500){\makebox(0,0){$\Theta_{3,1}$}}%
\put(50.6000,-10.5500){\makebox(0,0){$s_3$}}%
\end{picture}}%

Case  $\rm{(I)}$
 \end{center}
From the explicit formula for the height pairing, which is given in \cite{shioda90}, we have
\[
\langle P_0^\pm, P_0^\pm \rangle = \frac 13,\quad \langle P_1^\pm, P_1^\pm \rangle = \frac 13 ,\quad \langle P_2^\pm, P_2^\pm \rangle = \frac 13 ,
\]
\[
\langle P_3^\pm, P_3^\pm \rangle=\frac 12,\quad \langle P_1^\pm, P_2^\pm \rangle = \frac 16 \quad or -\frac 16.
\]
 If we denote $P_1:=P_1^+$, then we can choose $P_2$ from $P_2^+$ and $P_2^-$ such that $\langle P_1, P_2 \rangle= \dfrac{1}{6}$.
 Also we denote $P_3:=P_3^+$. Then $P_1$, $P_2$, and $P_3$ are a basis of the lattice $E(\CC(t))= A_1^*\oplus \dfrac{1}{6} \left [\begin{array}{cc}
                                                                2 & 1 \\
                                                                1 & 2
                                                                \end{array} \right ]$. We have $P_0\dot{+}P_1=P_2$ .
                                                                
Let $\mcC$ is weak contact conic, then $\mcC$ gives rise to section $s_{\mcC}^{\pm}$. We denote the corresponding 
element in $E_{\mcQ, z_0}(\CC(t))$ by $P_\mcC^\pm$.

\begin{lem}
If there exists a weak contact conic $\mcC$ of type $i$, $i=1\dots 6$, then $h_i=\langle P_{\mcC^\pm}, P_{\mcC^\pm} \rangle $ is given by the following table:
\setlength\extrarowheight{10pt}
\begin{center}
\begin{tabular}{|c|c|c|c|c|c|c|} \hline 
 & type 1 &type 2& type 3& type 4& type 5& type 6  \\ \hline  
$h_i $&$\dfrac{4}{3}$&$\dfrac{3}{2}$&$\dfrac{5}{6}$&$1$&$\dfrac{1}{3}$&$2$\\[10pt]\hline                                                               
\end{tabular}                                                                                             
\end{center}
\end{lem}                                                                
                                  
From the above table we can find weak contact conics as follows.
\begin{lem}
With the same notation as before, we have the following weak contact conics:
\begin{center}
\begin{tabular}{|c|c|c|} \hline 
type & corresponding elements of $E_{\mcQ,z_0}(\CC(t))$ & \#Weak CC\\ \hline
1 & $[2]P_1, [2]P_2, [2]P_0 $ & 3 \\[5pt] \hline 
2 & $P_0\dot{+}P_2\dot{+}P_3$, $P_0\dot{+}P_2\dot{+}(-P_3)$, & 4  \\
&  $P_1\dot{+}(-P_0)\dot{+}P_3$, $P_1\dot{+}(-P_0)\dot{+}(-P_3)$&  \\[5pt]  \hline
3 & $P_1\dot{+}P_3$,  $P_1\dot{+}(-P_3)$, $P_2\dot{+}P_3$, $P_2\dot{+}(-P_3)$, & 4 \\[5pt] \hline
4 & $P_1\dot{+}P_2$ & 1 \\[5pt] \hline
5 & $P_0$ & 1  \\[5pt] \hline
6 & $[2]P_3$ & 1 \\[5pt] \hline                                                               
\end{tabular}                                                                                             
\end{center}
\end{lem}

We will use this lemma in sections 4 and 5 to give examples of Zariski pairs.

{\bf Case $\mathrm{(II)}$}

 We label irreducible components of singular fibers of type $\I_2$ and $\I_3$ in the same way as in Case $\mathrm{(I)}$ and those of
  type $\I_3$ such that $\Theta_{\infty,1}$, $\Theta_{\infty,2}$ are irreducible 
 components from $l_{z_0}$.
 By our construction of $S_{\mcQ, z_o}$, we may assume that
 $s_i^+$ ($i = 1, 2$) meet each singular fiber as in the
 figure below.
  \begin{center}
{\unitlength 0.1in%
\begin{picture}( 44.9000, 21.8500)(  5.9000,-22.7000)%
%
\special{pn 13}%
\special{pa 1285 253}%
\special{pa 1256 273}%
\special{pa 1172 333}%
\special{pa 1145 353}%
\special{pa 1119 374}%
\special{pa 1093 396}%
\special{pa 1045 440}%
\special{pa 1023 463}%
\special{pa 1002 487}%
\special{pa 982 512}%
\special{pa 965 537}%
\special{pa 949 564}%
\special{pa 935 592}%
\special{pa 923 621}%
\special{pa 913 651}%
\special{pa 906 682}%
\special{pa 901 714}%
\special{pa 898 747}%
\special{pa 898 781}%
\special{pa 900 814}%
\special{pa 905 848}%
\special{pa 911 881}%
\special{pa 920 914}%
\special{pa 931 946}%
\special{pa 945 977}%
\special{pa 961 1006}%
\special{pa 979 1034}%
\special{pa 999 1060}%
\special{pa 1021 1083}%
\special{pa 1046 1105}%
\special{pa 1071 1125}%
\special{pa 1097 1145}%
\special{pa 1122 1164}%
\special{pa 1145 1185}%
\special{pa 1167 1206}%
\special{pa 1186 1229}%
\special{pa 1203 1254}%
\special{pa 1218 1279}%
\special{pa 1231 1306}%
\special{pa 1242 1334}%
\special{pa 1251 1363}%
\special{pa 1259 1393}%
\special{pa 1266 1425}%
\special{pa 1271 1456}%
\special{pa 1274 1489}%
\special{pa 1277 1523}%
\special{pa 1279 1557}%
\special{pa 1279 1627}%
\special{pa 1278 1662}%
\special{pa 1277 1698}%
\special{pa 1275 1734}%
\special{pa 1273 1771}%
\special{pa 1271 1807}%
\special{pa 1268 1844}%
\special{pa 1268 1846}%
\special{fp}%
%
\special{pn 13}%
\special{pa 1022 269}%
\special{pa 1044 296}%
\special{pa 1065 322}%
\special{pa 1087 349}%
\special{pa 1129 403}%
\special{pa 1169 457}%
\special{pa 1205 511}%
\special{pa 1223 539}%
\special{pa 1239 566}%
\special{pa 1254 594}%
\special{pa 1280 650}%
\special{pa 1291 678}%
\special{pa 1300 706}%
\special{pa 1308 735}%
\special{pa 1314 764}%
\special{pa 1318 793}%
\special{pa 1320 823}%
\special{pa 1321 852}%
\special{pa 1320 882}%
\special{pa 1318 912}%
\special{pa 1314 943}%
\special{pa 1310 973}%
\special{pa 1296 1035}%
\special{pa 1288 1066}%
\special{pa 1268 1128}%
\special{pa 1257 1159}%
\special{pa 1245 1191}%
\special{pa 1232 1222}%
\special{pa 1219 1254}%
\special{pa 1205 1286}%
\special{pa 1191 1317}%
\special{pa 1177 1349}%
\special{pa 1132 1445}%
\special{pa 1126 1457}%
\special{fp}%
%
\special{pn 13}%
\special{pa 1988 253}%
\special{pa 1959 273}%
\special{pa 1875 333}%
\special{pa 1848 353}%
\special{pa 1822 374}%
\special{pa 1796 396}%
\special{pa 1772 417}%
\special{pa 1748 440}%
\special{pa 1726 463}%
\special{pa 1705 487}%
\special{pa 1686 512}%
\special{pa 1668 537}%
\special{pa 1652 564}%
\special{pa 1638 592}%
\special{pa 1626 621}%
\special{pa 1616 651}%
\special{pa 1609 682}%
\special{pa 1604 714}%
\special{pa 1601 747}%
\special{pa 1601 781}%
\special{pa 1603 814}%
\special{pa 1607 848}%
\special{pa 1614 881}%
\special{pa 1623 914}%
\special{pa 1634 946}%
\special{pa 1648 977}%
\special{pa 1664 1006}%
\special{pa 1682 1034}%
\special{pa 1702 1060}%
\special{pa 1724 1083}%
\special{pa 1749 1105}%
\special{pa 1774 1125}%
\special{pa 1800 1144}%
\special{pa 1825 1164}%
\special{pa 1849 1184}%
\special{pa 1870 1206}%
\special{pa 1890 1228}%
\special{pa 1907 1253}%
\special{pa 1922 1278}%
\special{pa 1935 1305}%
\special{pa 1946 1333}%
\special{pa 1955 1362}%
\special{pa 1963 1392}%
\special{pa 1969 1424}%
\special{pa 1974 1455}%
\special{pa 1978 1488}%
\special{pa 1981 1522}%
\special{pa 1983 1590}%
\special{pa 1983 1626}%
\special{pa 1982 1661}%
\special{pa 1980 1697}%
\special{pa 1979 1734}%
\special{pa 1976 1770}%
\special{pa 1974 1807}%
\special{pa 1971 1843}%
\special{pa 1971 1846}%
\special{fp}%
%
\special{pn 13}%
\special{pa 1726 269}%
\special{pa 1748 296}%
\special{pa 1769 323}%
\special{pa 1791 349}%
\special{pa 1812 376}%
\special{pa 1832 403}%
\special{pa 1853 430}%
\special{pa 1891 484}%
\special{pa 1909 512}%
\special{pa 1926 539}%
\special{pa 1942 566}%
\special{pa 1957 594}%
\special{pa 1971 622}%
\special{pa 1983 650}%
\special{pa 1994 678}%
\special{pa 2003 707}%
\special{pa 2011 736}%
\special{pa 2017 765}%
\special{pa 2021 794}%
\special{pa 2023 823}%
\special{pa 2024 853}%
\special{pa 2023 883}%
\special{pa 2021 913}%
\special{pa 2017 943}%
\special{pa 2013 974}%
\special{pa 2006 1004}%
\special{pa 1999 1035}%
\special{pa 1991 1066}%
\special{pa 1981 1097}%
\special{pa 1971 1129}%
\special{pa 1960 1160}%
\special{pa 1948 1191}%
\special{pa 1922 1255}%
\special{pa 1909 1286}%
\special{pa 1894 1318}%
\special{pa 1880 1350}%
\special{pa 1835 1446}%
\special{pa 1830 1457}%
\special{fp}%
%
\special{pn 13}%
\special{pa 2387 253}%
\special{pa 2387 1851}%
\special{fp}%
%
\special{pn 13}%
\special{pa 4008 258}%
\special{pa 4014 290}%
\special{pa 4019 321}%
\special{pa 4031 385}%
\special{pa 4036 416}%
\special{pa 4046 480}%
\special{pa 4051 511}%
\special{pa 4059 575}%
\special{pa 4063 606}%
\special{pa 4069 670}%
\special{pa 4075 766}%
\special{pa 4075 862}%
\special{pa 4074 894}%
\special{pa 4066 1022}%
\special{pa 4063 1054}%
\special{pa 4061 1086}%
\special{pa 4052 1182}%
\special{pa 4050 1213}%
\special{pa 4047 1245}%
\special{pa 4043 1277}%
\special{pa 4028 1307}%
\special{pa 4002 1329}%
\special{pa 3971 1331}%
\special{pa 3938 1313}%
\special{pa 3911 1283}%
\special{pa 3897 1249}%
\special{pa 3901 1221}%
\special{pa 3926 1204}%
\special{pa 3962 1197}%
\special{pa 3996 1201}%
\special{pa 4023 1213}%
\special{pa 4044 1232}%
\special{pa 4059 1258}%
\special{pa 4069 1290}%
\special{pa 4076 1325}%
\special{pa 4078 1363}%
\special{pa 4079 1403}%
\special{pa 4078 1444}%
\special{pa 4076 1484}%
\special{pa 4074 1523}%
\special{pa 4072 1560}%
\special{pa 4070 1595}%
\special{pa 4069 1629}%
\special{pa 4068 1662}%
\special{pa 4066 1724}%
\special{pa 4064 1782}%
\special{pa 4064 1811}%
\special{pa 4063 1840}%
\special{pa 4063 1851}%
\special{fp}%
%
\special{pn 13}%
\special{pa 4397 264}%
\special{pa 4403 296}%
\special{pa 4408 327}%
\special{pa 4413 359}%
\special{pa 4419 391}%
\special{pa 4424 422}%
\special{pa 4434 486}%
\special{pa 4438 517}%
\special{pa 4443 549}%
\special{pa 4447 581}%
\special{pa 4450 613}%
\special{pa 4453 644}%
\special{pa 4456 676}%
\special{pa 4460 740}%
\special{pa 4462 804}%
\special{pa 4462 868}%
\special{pa 4460 932}%
\special{pa 4456 996}%
\special{pa 4454 1029}%
\special{pa 4451 1061}%
\special{pa 4449 1092}%
\special{pa 4440 1188}%
\special{pa 4438 1220}%
\special{pa 4435 1251}%
\special{pa 4430 1283}%
\special{pa 4415 1313}%
\special{pa 4389 1335}%
\special{pa 4357 1336}%
\special{pa 4325 1317}%
\special{pa 4299 1286}%
\special{pa 4285 1253}%
\special{pa 4291 1225}%
\special{pa 4317 1208}%
\special{pa 4353 1202}%
\special{pa 4387 1206}%
\special{pa 4414 1219}%
\special{pa 4434 1239}%
\special{pa 4448 1265}%
\special{pa 4458 1297}%
\special{pa 4464 1333}%
\special{pa 4466 1371}%
\special{pa 4466 1411}%
\special{pa 4465 1452}%
\special{pa 4462 1492}%
\special{pa 4460 1531}%
\special{pa 4459 1568}%
\special{pa 4457 1603}%
\special{pa 4456 1637}%
\special{pa 4454 1701}%
\special{pa 4454 1731}%
\special{pa 4453 1761}%
\special{pa 4453 1790}%
\special{pa 4452 1819}%
\special{pa 4452 1857}%
\special{fp}%
%
\special{pn 13}%
\special{pa 771 474}%
\special{pa 930 474}%
\special{fp}%
%
\special{pn 13}%
\special{pa 1072 474}%
\special{pa 1640 474}%
\special{fp}%
%
\special{pn 13}%
\special{pa 2447 685}%
\special{pa 3070 685}%
\special{fp}%
%
\special{pn 13}%
\special{pa 3160 250}%
\special{pa 3160 1805}%
\special{fp}%
%
\special{pn 13}%
\special{pa 2325 1355}%
\special{pa 2640 555}%
\special{fp}%
%
\special{pn 13}%
\special{pa 2655 625}%
\special{pa 2310 280}%
\special{fp}%
%
\special{pn 13}%
\special{pa 2430 475}%
\special{pa 3205 475}%
\special{fp}%
%
\special{pn 13}%
\special{pa 1785 475}%
\special{pa 2325 475}%
\special{fp}%
%
\special{pn 13}%
\special{pa 1665 860}%
\special{pa 2475 860}%
\special{fp}%
%
\special{pn 13}%
\special{pa 2565 860}%
\special{pa 3105 865}%
\special{fp}%
%
\special{pn 13}%
\special{pa 1650 685}%
\special{pa 2315 685}%
\special{fp}%
%
\special{pn 13}%
\special{pa 3080 1335}%
\special{pa 3545 950}%
\special{fp}%
%
\special{pn 13}%
\special{pa 3550 1000}%
\special{pa 3135 265}%
\special{fp}%
%
\special{pn 13}%
\special{pa 3290 470}%
\special{pa 4655 470}%
\special{fp}%
%
\special{pn 13}%
\special{pa 770 690}%
\special{pa 1230 680}%
\special{fp}%
%
\special{pn 13}%
\special{pa 1335 685}%
\special{pa 1555 680}%
\special{fp}%
%
\special{pn 13}%
\special{pa 3215 680}%
\special{pa 4655 685}%
\special{fp}%
%
\special{pn 13}%
\special{pa 760 860}%
\special{pa 855 865}%
\special{fp}%
%
\special{pn 13}%
\special{pa 950 865}%
\special{pa 1550 865}%
\special{fp}%
%
\special{pn 13}%
\special{pa 3220 865}%
\special{pa 4640 865}%
\special{fp}%
%
\special{pn 13}%
\special{pa 765 1045}%
\special{pa 940 1045}%
\special{fp}%
%
\special{pn 13}%
\special{pa 1020 1045}%
\special{pa 1955 1045}%
\special{fp}%
%
\special{pn 13}%
\special{pa 2050 1045}%
\special{pa 2350 1045}%
\special{fp}%
%
\special{pn 13}%
\special{pa 2405 1040}%
\special{pa 3100 1040}%
\special{fp}%
%
\special{pn 13}%
\special{pa 3230 1040}%
\special{pa 4640 1040}%
\special{fp}%
%
\special{pn 13}%
\special{pa 590 85}%
\special{pa 5080 85}%
\special{pa 5080 2270}%
\special{pa 590 2270}%
\special{pa 590 85}%
\special{pa 5080 85}%
\special{fp}%
\put(48.3500,-6.6500){\makebox(0,0){$s_1$}}%
\put(48.3500,-4.6000){\makebox(0,0){$s_0$}}%
\put(48.4000,-8.6500){\makebox(0,0){$s_2$}}%
\put(48.2500,-16.5000){\makebox(0,0){$O$}}%
\put(48.2500,-10.5000){\makebox(0,0){$s_3$}}%
%
\special{pn 13}%
\special{pa 750 1700}%
\special{pa 4650 1700}%
\special{fp}%
\put(14.5200,-18.0800){\makebox(0,0){$\Theta_{1,0}$}}%
\put(21.5200,-18.1800){\makebox(0,0){$\Theta_{2,0}$}}%
\put(33.4500,-18.2000){\makebox(0,0){$\Theta_{\infty,0}$}}%
\put(25.6200,-18.2300){\makebox(0,0){$\Theta_{3,0}$}}%
\put(14.0700,-12.1300){\makebox(0,0){$\Theta_{1,1}$}}%
\put(21.0700,-12.2300){\makebox(0,0){$\Theta_{2,1}$}}%
\put(34.8500,-11.6500){\makebox(0,0){$\Theta_{\infty,2}$}}%
\put(25.9400,-11.6800){\makebox(0,0){$\Theta_{3,2}$}}%
\put(25.7500,-3.8000){\makebox(0,0){$\Theta_{3,1}$}}%
\put(34.1000,-3.7500){\makebox(0,0){$\Theta_{\infty,1}$}}%
\end{picture}}%

  Case  $\mathrm{(II)}$
 \end{center}
 
From the explicit formula for the height pairing, we have
\[
\langle P_1^\pm, P_1^\pm \rangle = \frac 16 ,\quad \langle P_2^\pm, P_2^\pm \rangle = \frac 16 ,\quad \langle P_1^\pm, P_2^\pm \rangle = 0
\]

 If we denote $P_1:=P_1^+$, then we can choose $P_2$ from $P_2^+$ and $P_2^-$ such that $\langle P_1, P_2 \rangle= 0$. Then $P_1$ and $P_2$ are a basis of the lattice $E(\CC(t))=\langle 1/6 \rangle^{\oplus 2}, $. We have $P_0=P_1\dot{+}P_2, \quad P_3=P_1\dot{+}(-P_2)$ .
 Since every weak contact conic intersects with $\Theta_{\infty,0}$, we have the same height pairing table with the Case (I) for each type of weak contact conic. 
\begin{lem}
With the same notation as before, we have the following weak contact conics:
\begin{center}

\begin{tabular}{|c|c|c|c|} \hline 
type &  $h_i$ & corresponding elements of $E_{\mcQ,z_0}(\CC(t))$ &\# weak CC\\ \hline
1 & $\frac{4}{3}$ & $[2]P_1 \dot{+}[2]P_2 $ & 1\\\hline 
2 & $\frac{3}{2}$ & $[3]P_1$, $[3]P_2$ &  2 \\ \hline
3 & $\frac{5}{6}$ & $P_1\dot{+}[-2]P_2$,  $[2]P_1\dot{+}(-P_2)$ & 2 \\ \hline
4 & $1 $ & $\emptyset$ & 0 \\ \hline
5 & $\frac{1}{3}$ & $P_1\dot{+}P_2$ & 1  \\ \hline
6 & $2$ & $\emptyset$ & 0 \\ \hline                                                               
\end{tabular}                                                                                             
\end{center}
\end{lem}

{\bf Case  $\mathrm{(III)}$}

 Since $l_{z_0}$ goes through the cusp, there is no conic tangent to $z_0$ and passing through the cusp.  We label irreducible components of singular fibers of type $\I_2$ in the same way as in Case $\mathrm{(I)}$ and those of
  type $\I_5$ such that $\Theta_{\infty,i}$, $(i=1,2,3,4)$ are irreducible 
 components from $l_{z_0}$ and the cusp.
 We may assume that
 $s_i^+$ ($i = 1, 2$) meet each singular fiber as follows.
 
   \begin{center}
{\unitlength 0.1in%
\begin{picture}( 44.9000, 21.8500)(  5.9000,-22.7000)%
%
\special{pn 13}%
\special{pa 1285 253}%
\special{pa 1256 273}%
\special{pa 1172 333}%
\special{pa 1145 353}%
\special{pa 1119 374}%
\special{pa 1093 396}%
\special{pa 1045 440}%
\special{pa 1023 463}%
\special{pa 1002 487}%
\special{pa 982 512}%
\special{pa 965 537}%
\special{pa 949 564}%
\special{pa 935 592}%
\special{pa 923 621}%
\special{pa 913 651}%
\special{pa 906 682}%
\special{pa 901 714}%
\special{pa 898 747}%
\special{pa 898 781}%
\special{pa 900 814}%
\special{pa 905 848}%
\special{pa 911 881}%
\special{pa 920 914}%
\special{pa 931 946}%
\special{pa 945 977}%
\special{pa 961 1006}%
\special{pa 979 1034}%
\special{pa 999 1060}%
\special{pa 1021 1083}%
\special{pa 1046 1105}%
\special{pa 1071 1125}%
\special{pa 1097 1145}%
\special{pa 1122 1164}%
\special{pa 1145 1185}%
\special{pa 1167 1206}%
\special{pa 1186 1229}%
\special{pa 1203 1254}%
\special{pa 1218 1279}%
\special{pa 1231 1306}%
\special{pa 1242 1334}%
\special{pa 1251 1363}%
\special{pa 1259 1393}%
\special{pa 1266 1425}%
\special{pa 1271 1456}%
\special{pa 1274 1489}%
\special{pa 1277 1523}%
\special{pa 1279 1557}%
\special{pa 1279 1627}%
\special{pa 1278 1662}%
\special{pa 1277 1698}%
\special{pa 1275 1734}%
\special{pa 1273 1771}%
\special{pa 1271 1807}%
\special{pa 1268 1844}%
\special{pa 1268 1846}%
\special{fp}%
%
\special{pn 13}%
\special{pa 1022 269}%
\special{pa 1044 296}%
\special{pa 1065 322}%
\special{pa 1087 349}%
\special{pa 1129 403}%
\special{pa 1169 457}%
\special{pa 1205 511}%
\special{pa 1223 539}%
\special{pa 1239 566}%
\special{pa 1254 594}%
\special{pa 1280 650}%
\special{pa 1291 678}%
\special{pa 1300 706}%
\special{pa 1308 735}%
\special{pa 1314 764}%
\special{pa 1318 793}%
\special{pa 1320 823}%
\special{pa 1321 852}%
\special{pa 1320 882}%
\special{pa 1318 912}%
\special{pa 1314 943}%
\special{pa 1310 973}%
\special{pa 1296 1035}%
\special{pa 1288 1066}%
\special{pa 1268 1128}%
\special{pa 1257 1159}%
\special{pa 1245 1191}%
\special{pa 1232 1222}%
\special{pa 1219 1254}%
\special{pa 1205 1286}%
\special{pa 1191 1317}%
\special{pa 1177 1349}%
\special{pa 1132 1445}%
\special{pa 1126 1457}%
\special{fp}%
%
\special{pn 13}%
\special{pa 1988 253}%
\special{pa 1959 273}%
\special{pa 1875 333}%
\special{pa 1848 353}%
\special{pa 1822 374}%
\special{pa 1796 396}%
\special{pa 1772 417}%
\special{pa 1748 440}%
\special{pa 1726 463}%
\special{pa 1705 487}%
\special{pa 1686 512}%
\special{pa 1668 537}%
\special{pa 1652 564}%
\special{pa 1638 592}%
\special{pa 1626 621}%
\special{pa 1616 651}%
\special{pa 1609 682}%
\special{pa 1604 714}%
\special{pa 1601 747}%
\special{pa 1601 781}%
\special{pa 1603 814}%
\special{pa 1607 848}%
\special{pa 1614 881}%
\special{pa 1623 914}%
\special{pa 1634 946}%
\special{pa 1648 977}%
\special{pa 1664 1006}%
\special{pa 1682 1034}%
\special{pa 1702 1060}%
\special{pa 1724 1083}%
\special{pa 1749 1105}%
\special{pa 1774 1125}%
\special{pa 1800 1144}%
\special{pa 1825 1164}%
\special{pa 1849 1184}%
\special{pa 1870 1206}%
\special{pa 1890 1228}%
\special{pa 1907 1253}%
\special{pa 1922 1278}%
\special{pa 1935 1305}%
\special{pa 1946 1333}%
\special{pa 1955 1362}%
\special{pa 1963 1392}%
\special{pa 1969 1424}%
\special{pa 1974 1455}%
\special{pa 1978 1488}%
\special{pa 1981 1522}%
\special{pa 1983 1590}%
\special{pa 1983 1626}%
\special{pa 1982 1661}%
\special{pa 1980 1697}%
\special{pa 1979 1734}%
\special{pa 1976 1770}%
\special{pa 1974 1807}%
\special{pa 1971 1843}%
\special{pa 1971 1846}%
\special{fp}%
%
\special{pn 13}%
\special{pa 1726 269}%
\special{pa 1748 296}%
\special{pa 1769 323}%
\special{pa 1791 349}%
\special{pa 1812 376}%
\special{pa 1832 403}%
\special{pa 1853 430}%
\special{pa 1891 484}%
\special{pa 1909 512}%
\special{pa 1926 539}%
\special{pa 1942 566}%
\special{pa 1957 594}%
\special{pa 1971 622}%
\special{pa 1983 650}%
\special{pa 1994 678}%
\special{pa 2003 707}%
\special{pa 2011 736}%
\special{pa 2017 765}%
\special{pa 2021 794}%
\special{pa 2023 823}%
\special{pa 2024 853}%
\special{pa 2023 883}%
\special{pa 2021 913}%
\special{pa 2017 943}%
\special{pa 2013 974}%
\special{pa 2006 1004}%
\special{pa 1999 1035}%
\special{pa 1991 1066}%
\special{pa 1981 1097}%
\special{pa 1971 1129}%
\special{pa 1960 1160}%
\special{pa 1948 1191}%
\special{pa 1922 1255}%
\special{pa 1909 1286}%
\special{pa 1894 1318}%
\special{pa 1880 1350}%
\special{pa 1835 1446}%
\special{pa 1830 1457}%
\special{fp}%
%
\special{pn 13}%
\special{pa 4008 258}%
\special{pa 4014 290}%
\special{pa 4019 321}%
\special{pa 4031 385}%
\special{pa 4036 416}%
\special{pa 4046 480}%
\special{pa 4051 511}%
\special{pa 4059 575}%
\special{pa 4063 606}%
\special{pa 4069 670}%
\special{pa 4075 766}%
\special{pa 4075 862}%
\special{pa 4074 894}%
\special{pa 4066 1022}%
\special{pa 4063 1054}%
\special{pa 4061 1086}%
\special{pa 4052 1182}%
\special{pa 4050 1213}%
\special{pa 4047 1245}%
\special{pa 4043 1277}%
\special{pa 4028 1307}%
\special{pa 4002 1329}%
\special{pa 3971 1331}%
\special{pa 3938 1313}%
\special{pa 3911 1283}%
\special{pa 3897 1249}%
\special{pa 3901 1221}%
\special{pa 3926 1204}%
\special{pa 3962 1197}%
\special{pa 3996 1201}%
\special{pa 4023 1213}%
\special{pa 4044 1232}%
\special{pa 4059 1258}%
\special{pa 4069 1290}%
\special{pa 4076 1325}%
\special{pa 4078 1363}%
\special{pa 4079 1403}%
\special{pa 4078 1444}%
\special{pa 4076 1484}%
\special{pa 4074 1523}%
\special{pa 4072 1560}%
\special{pa 4070 1595}%
\special{pa 4069 1629}%
\special{pa 4068 1662}%
\special{pa 4066 1724}%
\special{pa 4064 1782}%
\special{pa 4064 1811}%
\special{pa 4063 1840}%
\special{pa 4063 1851}%
\special{fp}%
%
\special{pn 13}%
\special{pa 4397 264}%
\special{pa 4403 296}%
\special{pa 4408 327}%
\special{pa 4413 359}%
\special{pa 4419 391}%
\special{pa 4424 422}%
\special{pa 4434 486}%
\special{pa 4438 517}%
\special{pa 4443 549}%
\special{pa 4447 581}%
\special{pa 4450 613}%
\special{pa 4453 644}%
\special{pa 4456 676}%
\special{pa 4460 740}%
\special{pa 4462 804}%
\special{pa 4462 868}%
\special{pa 4460 932}%
\special{pa 4456 996}%
\special{pa 4454 1029}%
\special{pa 4451 1061}%
\special{pa 4449 1092}%
\special{pa 4440 1188}%
\special{pa 4438 1220}%
\special{pa 4435 1251}%
\special{pa 4430 1283}%
\special{pa 4415 1313}%
\special{pa 4389 1335}%
\special{pa 4357 1336}%
\special{pa 4325 1317}%
\special{pa 4299 1286}%
\special{pa 4285 1253}%
\special{pa 4291 1225}%
\special{pa 4317 1208}%
\special{pa 4353 1202}%
\special{pa 4387 1206}%
\special{pa 4414 1219}%
\special{pa 4434 1239}%
\special{pa 4448 1265}%
\special{pa 4458 1297}%
\special{pa 4464 1333}%
\special{pa 4466 1371}%
\special{pa 4466 1411}%
\special{pa 4465 1452}%
\special{pa 4462 1492}%
\special{pa 4460 1531}%
\special{pa 4459 1568}%
\special{pa 4457 1603}%
\special{pa 4456 1637}%
\special{pa 4454 1701}%
\special{pa 4454 1731}%
\special{pa 4453 1761}%
\special{pa 4453 1790}%
\special{pa 4452 1819}%
\special{pa 4452 1857}%
\special{fp}%
%
\special{pn 13}%
\special{pa 3290 470}%
\special{pa 4655 470}%
\special{fp}%
%
\special{pn 13}%
\special{pa 590 85}%
\special{pa 5080 85}%
\special{pa 5080 2270}%
\special{pa 590 2270}%
\special{pa 590 85}%
\special{pa 5080 85}%
\special{fp}%
%
\special{pn 13}%
\special{pa 760 860}%
\special{pa 855 865}%
\special{fp}%
%
\special{pn 13}%
\special{pa 3220 865}%
\special{pa 4640 865}%
\special{fp}%
%
\special{pn 13}%
\special{pa 3605 265}%
\special{pa 3611 297}%
\special{pa 3616 328}%
\special{pa 3628 392}%
\special{pa 3633 423}%
\special{pa 3643 487}%
\special{pa 3648 518}%
\special{pa 3656 582}%
\special{pa 3660 613}%
\special{pa 3666 677}%
\special{pa 3672 773}%
\special{pa 3672 869}%
\special{pa 3671 901}%
\special{pa 3663 1029}%
\special{pa 3660 1061}%
\special{pa 3658 1093}%
\special{pa 3649 1189}%
\special{pa 3647 1220}%
\special{pa 3644 1252}%
\special{pa 3640 1284}%
\special{pa 3625 1314}%
\special{pa 3599 1336}%
\special{pa 3568 1338}%
\special{pa 3535 1320}%
\special{pa 3508 1290}%
\special{pa 3494 1256}%
\special{pa 3498 1228}%
\special{pa 3523 1211}%
\special{pa 3559 1204}%
\special{pa 3593 1208}%
\special{pa 3620 1220}%
\special{pa 3641 1239}%
\special{pa 3656 1265}%
\special{pa 3666 1297}%
\special{pa 3673 1332}%
\special{pa 3675 1370}%
\special{pa 3676 1410}%
\special{pa 3675 1451}%
\special{pa 3673 1491}%
\special{pa 3671 1530}%
\special{pa 3669 1567}%
\special{pa 3667 1602}%
\special{pa 3666 1636}%
\special{pa 3665 1669}%
\special{pa 3663 1731}%
\special{pa 3661 1789}%
\special{pa 3661 1818}%
\special{pa 3660 1847}%
\special{pa 3660 1858}%
\special{fp}%
%
\special{pn 13}%
\special{pa 2520 285}%
\special{pa 2520 1780}%
\special{fp}%
%
\special{pn 13}%
\special{pa 775 505}%
\special{pa 934 505}%
\special{fp}%
%
\special{pn 13}%
\special{pa 1050 505}%
\special{pa 1618 505}%
\special{fp}%
%
\special{pn 13}%
\special{pa 2440 1205}%
\special{pa 3160 1280}%
\special{fp}%
\special{pa 2965 1335}%
\special{pa 3380 760}%
\special{fp}%
\special{pa 3375 835}%
\special{pa 2995 510}%
\special{fp}%
%
\special{pn 13}%
\special{pa 1745 505}%
\special{pa 2465 505}%
\special{fp}%
\special{pa 2575 505}%
\special{pa 3290 465}%
\special{fp}%
%
\special{pn 13}%
\special{pa 2435 400}%
\special{pa 3225 645}%
\special{fp}%
%
\special{pn 13}%
\special{pa 960 865}%
\special{pa 1960 865}%
\special{fp}%
\special{pa 2070 865}%
\special{pa 2450 865}%
\special{fp}%
%
\special{pn 13}%
\special{pa 2565 865}%
\special{pa 3230 860}%
\special{fp}%
%
\special{pn 13}%
\special{pa 765 1105}%
\special{pa 1225 1105}%
\special{fp}%
\special{pa 1325 1105}%
\special{pa 1705 1105}%
\special{fp}%
\special{pa 1805 1105}%
\special{pa 2460 1105}%
\special{fp}%
\special{pa 2560 1105}%
\special{pa 4645 1105}%
\special{fp}%
%
\special{pn 13}%
\special{pa 750 1600}%
\special{pa 4690 1600}%
\special{fp}%
\put(47.9000,-4.4500){\makebox(0,0){$s_1$}}%
\put(47.8500,-8.3500){\makebox(0,0){$s_2$}}%
\put(47.7500,-10.7500){\makebox(0,0){$s_3$}}%
\put(14.2500,-18.2000){\makebox(0,0){$\Theta_{1,0}$}}%
\put(21.2000,-18.1500){\makebox(0,0){$\Theta_{2,0}$}}%
\put(27.1000,-18.2500){\makebox(0,0){$\Theta_{\infty,0}$}}%
\put(47.9500,-16.0000){\makebox(0,0){$O$}}%
\put(14.3000,-12.2500){\makebox(0,0){$\Theta_{1,1}$}}%
\put(21.1500,-12.4500){\makebox(0,0){$\Theta_{2,1}$}}%
\put(27.6500,-3.7500){\makebox(0,0){$\Theta_{\infty,1}$}}%
\put(34.1500,-6.7500){\makebox(0,0){$\Theta_{\infty,2}$}}%
\put(33.7000,-10.1500){\makebox(0,0){$\Theta_{\infty,3}$}}%
\put(27.6000,-13.7000){\makebox(0,0){$\Theta_{\infty,4}$}}%
\end{picture}}%

  Case  $\mathrm{(III)}$
 \end{center}
 
Hence we have
\[
\langle P_1^\pm, P_1^\pm \rangle = \frac 15 ,\quad \langle P_2^\pm, P_2^\pm \rangle = \frac{3}{10} ,\quad \langle P_3^\pm, P_3^\pm \rangle = \frac{3}{10} , \]

\[
\langle P_1^\pm, P_2^\pm \rangle = \frac{1}{10} \quad or \quad -\frac{1}{10}
\]
\setlength\extrarowheight{0pt}
 If we denote $P_1:=P_1^+$, we can choose $P_2$ from $P_2^+$ and $P_2^-$ such that $\langle P_1, P_2 \rangle= \dfrac{1}{10}$. Then    $P_1$ and $P_2$ are a basis of the lattice $E(\CC(t))=\dfrac{1}{10} \left [ \begin{array}{cc}   
                                                                2 & 1 \\
                                                                1 & 3
                                                                \end{array} \right ] $. We have $P_3=P_1\dot{+}(-P_2) $ .
Likewise previous cases, we have the same height pairing table. 
\begin{lem}
With the same notation as before, we have the following weak contact conics:
\begin{center}
\begin{tabular}{|c|c|c|c|} \hline 
type & $h_i$ & corresponding elements of $E_{\mcQ,z_0}(\CC(t))$ &\# weak CC\\ \hline
2 & $\frac{3}{2}$ & $[2]P_1\dot{+}P_2$, $P_2\dot{+}[-3]P_1$ &  2 \\[5pt]  \hline
4 & 1 & $P_1\dot{+}[-2]P_2$ & 1 \\[5pt] \hline
6 & 2 & $\emptyset$ & 0 \\[5pt] \hline                                        
\end{tabular}                                                                                             
\end{center}
\end{lem}
Note that since no conic tangent to $z_0$ could not pass through the cusp, we omit types 1,3 and 5.

{\bf Case  $\mathrm{(IV)}$}

Let $l_{z_0}$ goes through a node $x_1$ and line $\mcL_1$ goes through $x_1$ and the cusp.
We label irreducible components of singular fibers of type $\I_2$ and $\I_3$ in the same way as before and those of
  type $\I_4$ such that $\Theta_{\infty,i}$, $(i=1,2,3)$ are irreducible 
 components from $l_{z_0}$ and $x_1$.
 We may assume that
 $s_i^+$ ($i = 1, 2$) meet each singular fiber as follows.
 
  \begin{center}
{\unitlength 0.1in%
\begin{picture}( 44.9000, 21.8500)(  5.9000,-22.7000)%
%
\special{pn 13}%
\special{pa 1285 253}%
\special{pa 1256 273}%
\special{pa 1172 333}%
\special{pa 1145 353}%
\special{pa 1119 374}%
\special{pa 1093 396}%
\special{pa 1045 440}%
\special{pa 1023 463}%
\special{pa 1002 487}%
\special{pa 982 512}%
\special{pa 965 537}%
\special{pa 949 564}%
\special{pa 935 592}%
\special{pa 923 621}%
\special{pa 913 651}%
\special{pa 906 682}%
\special{pa 901 714}%
\special{pa 898 747}%
\special{pa 898 781}%
\special{pa 900 814}%
\special{pa 905 848}%
\special{pa 911 881}%
\special{pa 920 914}%
\special{pa 931 946}%
\special{pa 945 977}%
\special{pa 961 1006}%
\special{pa 979 1034}%
\special{pa 999 1060}%
\special{pa 1021 1083}%
\special{pa 1046 1105}%
\special{pa 1071 1125}%
\special{pa 1097 1145}%
\special{pa 1122 1164}%
\special{pa 1145 1185}%
\special{pa 1167 1206}%
\special{pa 1186 1229}%
\special{pa 1203 1254}%
\special{pa 1218 1279}%
\special{pa 1231 1306}%
\special{pa 1242 1334}%
\special{pa 1251 1363}%
\special{pa 1259 1393}%
\special{pa 1266 1425}%
\special{pa 1271 1456}%
\special{pa 1274 1489}%
\special{pa 1277 1523}%
\special{pa 1279 1557}%
\special{pa 1279 1627}%
\special{pa 1278 1662}%
\special{pa 1277 1698}%
\special{pa 1275 1734}%
\special{pa 1273 1771}%
\special{pa 1271 1807}%
\special{pa 1268 1844}%
\special{pa 1268 1846}%
\special{fp}%
%
\special{pn 13}%
\special{pa 1022 269}%
\special{pa 1044 296}%
\special{pa 1065 322}%
\special{pa 1087 349}%
\special{pa 1129 403}%
\special{pa 1169 457}%
\special{pa 1205 511}%
\special{pa 1223 539}%
\special{pa 1239 566}%
\special{pa 1254 594}%
\special{pa 1280 650}%
\special{pa 1291 678}%
\special{pa 1300 706}%
\special{pa 1308 735}%
\special{pa 1314 764}%
\special{pa 1318 793}%
\special{pa 1320 823}%
\special{pa 1321 852}%
\special{pa 1320 882}%
\special{pa 1318 912}%
\special{pa 1314 943}%
\special{pa 1310 973}%
\special{pa 1296 1035}%
\special{pa 1288 1066}%
\special{pa 1268 1128}%
\special{pa 1257 1159}%
\special{pa 1245 1191}%
\special{pa 1232 1222}%
\special{pa 1219 1254}%
\special{pa 1205 1286}%
\special{pa 1191 1317}%
\special{pa 1177 1349}%
\special{pa 1132 1445}%
\special{pa 1126 1457}%
\special{fp}%
%
\special{pn 13}%
\special{pa 4008 258}%
\special{pa 4014 290}%
\special{pa 4019 321}%
\special{pa 4031 385}%
\special{pa 4036 416}%
\special{pa 4046 480}%
\special{pa 4051 511}%
\special{pa 4059 575}%
\special{pa 4063 606}%
\special{pa 4069 670}%
\special{pa 4075 766}%
\special{pa 4075 862}%
\special{pa 4074 894}%
\special{pa 4066 1022}%
\special{pa 4063 1054}%
\special{pa 4061 1086}%
\special{pa 4052 1182}%
\special{pa 4050 1213}%
\special{pa 4047 1245}%
\special{pa 4043 1277}%
\special{pa 4028 1307}%
\special{pa 4002 1329}%
\special{pa 3971 1331}%
\special{pa 3938 1313}%
\special{pa 3911 1283}%
\special{pa 3897 1249}%
\special{pa 3901 1221}%
\special{pa 3926 1204}%
\special{pa 3962 1197}%
\special{pa 3996 1201}%
\special{pa 4023 1213}%
\special{pa 4044 1232}%
\special{pa 4059 1258}%
\special{pa 4069 1290}%
\special{pa 4076 1325}%
\special{pa 4078 1363}%
\special{pa 4079 1403}%
\special{pa 4078 1444}%
\special{pa 4076 1484}%
\special{pa 4074 1523}%
\special{pa 4072 1560}%
\special{pa 4070 1595}%
\special{pa 4069 1629}%
\special{pa 4068 1662}%
\special{pa 4066 1724}%
\special{pa 4064 1782}%
\special{pa 4064 1811}%
\special{pa 4063 1840}%
\special{pa 4063 1851}%
\special{fp}%
%
\special{pn 13}%
\special{pa 4397 264}%
\special{pa 4403 296}%
\special{pa 4408 327}%
\special{pa 4413 359}%
\special{pa 4419 391}%
\special{pa 4424 422}%
\special{pa 4434 486}%
\special{pa 4438 517}%
\special{pa 4443 549}%
\special{pa 4447 581}%
\special{pa 4450 613}%
\special{pa 4453 644}%
\special{pa 4456 676}%
\special{pa 4460 740}%
\special{pa 4462 804}%
\special{pa 4462 868}%
\special{pa 4460 932}%
\special{pa 4456 996}%
\special{pa 4454 1029}%
\special{pa 4451 1061}%
\special{pa 4449 1092}%
\special{pa 4440 1188}%
\special{pa 4438 1220}%
\special{pa 4435 1251}%
\special{pa 4430 1283}%
\special{pa 4415 1313}%
\special{pa 4389 1335}%
\special{pa 4357 1336}%
\special{pa 4325 1317}%
\special{pa 4299 1286}%
\special{pa 4285 1253}%
\special{pa 4291 1225}%
\special{pa 4317 1208}%
\special{pa 4353 1202}%
\special{pa 4387 1206}%
\special{pa 4414 1219}%
\special{pa 4434 1239}%
\special{pa 4448 1265}%
\special{pa 4458 1297}%
\special{pa 4464 1333}%
\special{pa 4466 1371}%
\special{pa 4466 1411}%
\special{pa 4465 1452}%
\special{pa 4462 1492}%
\special{pa 4460 1531}%
\special{pa 4459 1568}%
\special{pa 4457 1603}%
\special{pa 4456 1637}%
\special{pa 4454 1701}%
\special{pa 4454 1731}%
\special{pa 4453 1761}%
\special{pa 4453 1790}%
\special{pa 4452 1819}%
\special{pa 4452 1857}%
\special{fp}%
%
\special{pn 13}%
\special{pa 3290 470}%
\special{pa 4655 470}%
\special{fp}%
%
\special{pn 13}%
\special{pa 590 85}%
\special{pa 5080 85}%
\special{pa 5080 2270}%
\special{pa 590 2270}%
\special{pa 590 85}%
\special{pa 5080 85}%
\special{fp}%
%
\special{pn 13}%
\special{pa 3220 865}%
\special{pa 4640 865}%
\special{fp}%
%
\special{pn 13}%
\special{pa 3605 265}%
\special{pa 3611 297}%
\special{pa 3616 328}%
\special{pa 3628 392}%
\special{pa 3633 423}%
\special{pa 3643 487}%
\special{pa 3648 518}%
\special{pa 3656 582}%
\special{pa 3660 613}%
\special{pa 3666 677}%
\special{pa 3672 773}%
\special{pa 3672 869}%
\special{pa 3671 901}%
\special{pa 3663 1029}%
\special{pa 3660 1061}%
\special{pa 3658 1093}%
\special{pa 3649 1189}%
\special{pa 3647 1220}%
\special{pa 3644 1252}%
\special{pa 3640 1284}%
\special{pa 3625 1314}%
\special{pa 3599 1336}%
\special{pa 3568 1338}%
\special{pa 3535 1320}%
\special{pa 3508 1290}%
\special{pa 3494 1256}%
\special{pa 3498 1228}%
\special{pa 3523 1211}%
\special{pa 3559 1204}%
\special{pa 3593 1208}%
\special{pa 3620 1220}%
\special{pa 3641 1239}%
\special{pa 3656 1265}%
\special{pa 3666 1297}%
\special{pa 3673 1332}%
\special{pa 3675 1370}%
\special{pa 3676 1410}%
\special{pa 3675 1451}%
\special{pa 3673 1491}%
\special{pa 3671 1530}%
\special{pa 3669 1567}%
\special{pa 3667 1602}%
\special{pa 3666 1636}%
\special{pa 3665 1669}%
\special{pa 3663 1731}%
\special{pa 3661 1789}%
\special{pa 3661 1818}%
\special{pa 3660 1847}%
\special{pa 3660 1858}%
\special{fp}%
%
\special{pn 13}%
\special{pa 2520 285}%
\special{pa 2520 1780}%
\special{fp}%
%
\special{pn 13}%
\special{pa 775 505}%
\special{pa 934 505}%
\special{fp}%
%
\special{pn 13}%
\special{pa 1050 505}%
\special{pa 1618 505}%
\special{fp}%
%
\special{pn 13}%
\special{pa 1745 505}%
\special{pa 2465 505}%
\special{fp}%
\special{pa 2575 505}%
\special{pa 3290 465}%
\special{fp}%
%
\special{pn 13}%
\special{pa 2565 865}%
\special{pa 3230 860}%
\special{fp}%
%
\special{pn 13}%
\special{pa 750 1600}%
\special{pa 4690 1600}%
\special{fp}%
\put(47.9000,-4.4500){\makebox(0,0){$s_2$}}%
\put(47.8500,-8.3500){\makebox(0,0){$s_1$}}%
\put(47.7500,-10.7500){\makebox(0,0){$s_3$}}%
\put(14.2500,-18.2000){\makebox(0,0){$\Theta_{1,0}$}}%
\put(21.2000,-18.1500){\makebox(0,0){$\Theta_{2,0}$}}%
\put(27.1000,-18.2500){\makebox(0,0){$\Theta_{\infty,0}$}}%
\put(47.9500,-16.0000){\makebox(0,0){$O$}}%
\put(14.3000,-12.2500){\makebox(0,0){$\Theta_{1,1}$}}%
\put(20.6000,-5.9000){\makebox(0,0){$\Theta_{2,1}$}}%
\put(27.6500,-3.7500){\makebox(0,0){$\Theta_{\infty,1}$}}%
\put(32.5500,-9.6000){\makebox(0,0){$\Theta_{\infty,2}$}}%
\put(27.8500,-12.8000){\makebox(0,0){$\Theta_{\infty,3}$}}%
%
\special{pn 13}%
\special{pa 2430 340}%
\special{pa 3130 775}%
\special{fp}%
\special{pa 3025 650}%
\special{pa 3035 1295}%
\special{fp}%
\special{pa 2315 1225}%
\special{pa 3135 1180}%
\special{fp}%
%
\special{pn 13}%
\special{pa 1680 245}%
\special{pa 1680 1840}%
\special{fp}%
%
\special{pn 13}%
\special{pa 2150 980}%
\special{pa 1635 275}%
\special{fp}%
\special{pa 1550 1310}%
\special{pa 2165 910}%
\special{fp}%
%
\special{pn 13}%
\special{pa 765 865}%
\special{pa 1245 865}%
\special{fp}%
\special{pa 1375 865}%
\special{pa 1620 865}%
\special{fp}%
\special{pa 1730 865}%
\special{pa 2470 865}%
\special{fp}%
%
\special{pn 13}%
\special{pa 760 1105}%
\special{pa 985 1105}%
\special{fp}%
\special{pa 1095 1105}%
\special{pa 1805 1105}%
\special{fp}%
\special{pa 1915 1105}%
\special{pa 2465 1105}%
\special{fp}%
\special{pa 2565 1105}%
\special{pa 4650 1105}%
\special{fp}%
\put(19.5500,-12.3500){\makebox(0,0){$\Theta_{2,2}$}}%
\end{picture}}%

  Case  $\mathrm{(IV)}$
 \end{center} 
Hence we have
\[
\langle P_1^\pm, P_1^\pm \rangle = \frac 13 ,~ \langle P_2^\pm, P_2^\pm \rangle = \frac{1}{12} ,~ \langle P_3^\pm, P_3^\pm \rangle = \frac{1}{2} ,~ \langle P_2^\pm, P_3^\pm \rangle = 0
\] 
 If we denote $P_1:=P_1^+$, $P_2:=P_2^+$ and $P_3:=P_3^+$, then $P_2$ and $P_3$ are a basis of the lattice $E(\CC(t))=  A_1^* \oplus \langle 1/12 \rangle $. We have $P_1=[-2]P_2 $ .

\begin{lem}
With the same notation as before, we have the following weak contact conics:
\begin{center}
\begin{tabular}{|c|c|c|c|} \hline 
type & $h_i$ & corresponding elements of $E_{\mcQ,z_0}(\CC(t))$ & \# weak CC\\ \hline
1 & $\frac{4}{3}$ & $[4]P_2 $ & 1\\ \hline
2 & $\frac{3}{2}$ & $\emptyset $ &  0 \\  \hline
3 & $\frac{5}{6}$ & $P_3\dot{+}[2]P_2$,  $P_3\dot{+}[-2]P_2$ & 2 \\ \hline
6 & 2 & $[2]P_3$ & 1 \\ \hline                                                               
\end{tabular}                                                                                             
\end{center}
\end{lem}
Note that since no conic tangent to $z_0$ could not pass through $2$ nodes, we omit type 4 and 5.

\section{Zariski pair}
Let $\mcQ$ and $\mcL_j$, $j=1,2,3$, be the quartic and lines as we defined in the previous section. For $i=0,1,2$, let $\mcC_i$ be the conics that correspond the rational points $[2]P_i$ $(i=0,1,2)$ of the Lemma 3.2, respectively. In particular, $[2]P_1$ and $[2]P_2$ arise from the line $\mcL_1$ and $\mcL_2$, respectively and $[2]P_0$ arises from the weak contact conic $\bar{\mcC}$ that goes through the 2 nodes and the cusp of $\mcQ$.
The following corollary imply from the Theorem~\ref{thm:homeo} and \ref{thm:homeo1}.

\begin{cor}\label{cor:ZariskiPair}There exists no homeomorphism between $(\PP^2,B_{1})$ and $(\PP^2,B_{2})$, if $(B_1,B_2)$ is one of the following:
\begin{eqnarray*}
(\mcB_{1,1},\mcB_{2,1})&=&(\mcQ+\mcL_1+\mcC_1,\mcQ+\mcL_2+\mcC_1)\\ 
(\mcB_{2,2},\mcB_{1,2})&=&(\mcQ+\mcL_2+\mcC_2,\mcQ+\mcL_1+\mcC_2)\\ 
\end{eqnarray*}

\begin{eqnarray*}
(\mcB_{1,1},\mcB_{1,2})&=&(\mcQ+\mcL_1+\mcC_1,\mcQ+\mcL_1+\mcC_2)\\ 
(\mcB_{1,1},\mcB_{1,0})&=&(\mcQ+\mcL_1+\mcC_1,\mcQ+\mcL_1+\mcC_0)\\ 
(\mcB_{2,2},\mcB_{2,1})&=&(\mcQ+\mcL_2+\mcC_2,\mcQ+\mcL_2+\mcC_1)\\
(\mcB_{2,2},\mcB_{2,0})&=&(\mcQ+\mcL_2+\mcC_2,\mcQ+\mcL_2+\mcC_0)\\
(\mcD_0,\mcD_1)&=&(\mcQ+\bar{\mcC}+\mcC_0,\mcQ+\bar{\mcC}+\mcC_1)\\
(\mcD_0,\mcD_2)&=&(\mcQ+\bar{\mcC}+\mcC_0,\mcQ+\bar{\mcC}+\mcC_2).
\end{eqnarray*}
\end{cor}
Proof.
 For the pairs $(\mcB_{1,1},\mcB_{2,1})$, we have $\mcB_{1,1}=\mcQ+\mcL_1+\mcC_1$ and $\mcB_{1,2}=\mcQ+\mcL_1+\mcC_2$. From section 3, we have $\mcL_1=f_{\mcQ,z_o}(s_1)$, $\mcC_1=f_{\mcQ,z_o}([2]s_1)$ and $\mcC_2=f_{\mcQ,z_o}([2]s_2)$. By Theorem~\ref{thm:homeo1}, there exists no homeomorphism between $(\PP^2,\mcB_{1,1})$ and $(\PP^2,\mcB_{2,1})$.
 For the remaining cases, it can be proved similarly. 
 
 In the next section, we will give explicit examples of Zariski pair.

\section{Example}
 We end this paper by giving explicit examples for an irreducible quartics with 2 nodes and 1 cusp and their contact conics observed so far,
 by which we have some examples of Zariski pairs.
 As for homogeneous coordinates of $\mathbb{P}^2$ we keep our previous notation, $[T,X,Z]$.
 
Let $C$ be a conic given by the equation $XZ-T^2=0$. Let $Q$ denote the standard quadratic transformations with respect to three lines given by the following equations: $T-X+Z=0, T+X-Z=0$, $Z=0$.

Let $L: X=0$ be tangent line at $p=[0,1,0]$ to $C$. Let us denote $\mcQ:=Q(C)$, $\mcC:=Q(L)$ and $z_0:=Q(L)$. Note that $\mcQ$ is quartic and $\mcC$ is conic. Let $l_{z_0}$ be the tangent line to $\mcQ$ at $z_0$. Then we have the equations of $\mcQ$, $\mcC$, $z_0$ and $l_{z_0}$ as follows:
\begin{eqnarray*}
F_{\mcQ} &=& Z^2X^2+2Z^2XT+Z^2T^2+2TX^2Z-2T^2XZ-4T^2X^2  \\
F_{\mcC} &=& 2TX+TZ-XZ\\
z_0 &=& [-1,1,-1],\\
F_{l_{z_0}} &=& T-X-2Z.
\end{eqnarray*}

We see that, $\mcC$ is conic tangent to $\mcQ$ at $z_0$ and passing through 2 nodes and cusp, and $l_{z_0}$ meets $\mcQ$ with two  distinct points. Let $\Phi$ be a coordinate change such that $l_{z_0}$ is
transformed into the line $Z=0$ and $z_0$ is mapped to $[0,1,0]$. Then
 $\Phi(\mcQ)$ and $\Phi(\mcC)$ are given by the affine equations as follows:
\begin{eqnarray*}
F_{\Phi(\mcQ)} &=& x^3+\dfrac{2t^2-3t}{2}x^2+(t^2-t^3)x+\dfrac{1}{8}t^2(t-1)^2=0 \\
F_{\Phi(\mcC)} &=& x+\dfrac{t^2}{2}-\dfrac{t}{2}=0,
\end{eqnarray*}
where $t=T/Z$ and $x=X/Z$.

Note that $\Phi(\mcQ)$ has 2-nodes at $[-1,-1,1], [1,0,1]$ and cusp at $[0,0,1]$. As we have discussed in section 3, we can assume line passing through 2 nodes and 2 lines passing through cusp and one of nodes  are basis of $ E(\mathbb{C}(t))$.
More precisely we have the table below: 
\begin{center}
 \begin{tabular}{l l}
 Equations & Element of $E_{\mcQ,z_0}(\CC(t))$ \\
 \hline

 $x=0$ & $P_1=\left(0,\dfrac{\sqrt{2}t(t-1)}{4}\right)$\\[12pt]

 $x-t=0$ & $P_2=\left(t,\dfrac{\sqrt{2} (t+1)t}{4}\right)$\\[12pt]
 
 $2x-t+1=0$ & $P_3=\left(\dfrac{t-1}{2},\dfrac{\sqrt{-2} (t-1)(t+1)}{4}\right)$\\[12pt]
 
 $2x+t^2-t=0$ & $P_0=\left( \dfrac{t-t^2}{2},\dfrac{\sqrt{2}t(t-1)(t+1)}{4} \right)$ \\
 
 \end{tabular}
\end{center}

Also we have $P_0=P_2\dot{+}(-P_1)$.

Next we will give explicit examples of weak contact conics as stated in section 3, Case(I)

{\rm type 1}. We have 3 weak contact conics passing through cusp only and tangent to $z_0$: $[2]P_1$, $[2]P_2,$ and $[2]P_0$.
\begin{center}
 \begin{tabular}{l l}
  Element of $E_{\mcQ,z_0}(\CC(t))$ & Conics\\
 \hline
 $[2]P_1^+=\left(\dfrac{(2t+3)t}{2}, \dfrac{\sqrt{2}(4t^2+5t+1)t}{4}\right)$ & $2x-(2t+3)t=0$\\[12pt]
 
 $[2]P_2^+=\left(\dfrac{(2t-1)t}{2}, -\dfrac{\sqrt{2}(4t^2-5t+1)t}{4}\right)$ & $2x-(2t-1)t=0$ \\[12pt]
 
 $[2]P_0^+=\left( \dfrac{t(t+4)}{8},-\dfrac{\sqrt{2}t(3t^2-8)}{32}\right)$ & $8x-t(t+4)=0$\\
 
 \end{tabular}
\end{center}

Let $\mcQ,\mcL_1,\mcL_2,\mcC_1,\mcC_2, \mcC_0, \bar{\mcC}$ are given as follows:
\begin{eqnarray*}
\mcQ: & & x^3+\dfrac{2t^2-3t}{2}x^2+(t^2-t^3)x+\dfrac{1}{8}t^2(t-1)^2=0\\
\mcL_1: & & x=0 \\
\mcL_2: & & x-t=0 \\
\mcC_1: & & 2x-(2t+3)t=0 \\
\mcC_2: & & 2x-(2t-1)t=0 \\
\mcC_0: & & 8x-t(t+4)=0\\
\bar{\mcC}: & & 2x+t^2-t=0 
\end{eqnarray*}
The fact that each of the following pairs has the same combinatorics, together with Corollary~\ref{cor:ZariskiPair} provide that the followings are Zariski pairs:
\begin{eqnarray*}
(\mcB_{1,1},\mcB_{2,1})&=&(\mcQ+\mcL_1+\mcC_1,\mcQ+\mcL_2+\mcC_1)\\ 
(\mcB_{2,2},\mcB_{1,2})&=&(\mcQ+\mcL_2+\mcC_2,\mcQ+\mcL_1+\mcC_2)\\ 
(\mcB_{1,1},\mcB_{1,2})&=&(\mcQ+\mcL_1+\mcC_1,\mcQ+\mcL_1+\mcC_2)\\ 
(\mcB_{1,1},\mcB_{1,0})&=&(\mcQ+\mcL_1+\mcC_1,\mcQ+\mcL_1+\mcC_0)\\ 
\end{eqnarray*}
\begin{eqnarray*}
(\mcB_{2,2},\mcB_{2,1})&=&(\mcQ+\mcL_2+\mcC_2,\mcQ+\mcL_2+\mcC_1)\\
(\mcB_{2,2},\mcB_{2,0})&=&(\mcQ+\mcL_2+\mcC_2,\mcQ+\mcL_2+\mcC_0)\\
(\mcD_0,\mcD_1)&=&(\mcQ+\bar{\mcC}+\mcC_0,\mcQ+\bar{\mcC}+\mcC_1)\\
(\mcD_0,\mcD_2)&=&(\mcQ+\bar{\mcC}+\mcC_0,\mcQ+\bar{\mcC}+\mcC_2).
\end{eqnarray*}

\noindent Khulan Tumenbayar\\%
Department of Mathematics\\%
School of Arts and Sciences,\\%
National University of Mongolia\\%
Ikh Surguuliin gudamj-1,Ulaanbaatar, MONGOLIA \\%
{\tt tkhulan@num.edu.mn}

\end{document}